\documentclass[twocolumn,amsthm,final]{autart}    

\usepackage{amsmath, amssymb, bbm, xspace}
\usepackage{amsfonts,mathrsfs}
\usepackage{mathtools}
\usepackage{color}
\usepackage{graphicx}
\usepackage{epsfig}
\usepackage{algorithm}
\usepackage{algorithmicx}
\usepackage[acronym]{glossaries}
\usepackage{subfigure}
\usepackage{cancel}
\usepackage{empheq}
\usepackage{placeins}
\usepackage{enumitem}
\usepackage{textgreek}
\usepackage{natbib}

\newcommand*{\QED}{\null\nobreak\hfill\ensuremath{\blacksquare}}%
\newcommand*{\QEDopen}{\null\nobreak\hfill\ensuremath{\square}}%
\newcommand{\bs}{\boldsymbol}
\newcommand{\mc}{\mathcal}
\renewcommand{\emph}{\textit}

\makeglossaries
\newacronym{GNEP}{GNEP}{generalized Nash equilibrium problem}
\newacronym{NE}{NE}{Nash equilibrium}
\newacronym{NEP}{NEP}{Nash equilibrium problem}
\newacronym{GNE}{GNE}{generalized Nash equilibrium}
\newacronym{v-GNE}{v-GNE}{variational \gls{GNE}}
\newacronym{ISS}{ISS}{input-to-state-stable}
\newacronym{PPA}{PPA}{proximal-point algorithm}
\newacronym{PPPA}{PPPA}{preconditioned \gls{PPA}}
\newacronym{VI}{VI}{variational inequality}
\newacronym{GAE}{GAE}{generalized aggregative equilibrium}
\newacronym{v-GAE}{v-GAE}{variational \gls{GAE}}
\newacronym{KKT}{KKT}{Karush--Kuhn--Tucker}
\newacronym{FQNE}{FQNE}{firmly quasinonexpansive}
\newacronym{FNE}{FNE}{firmly nonexpansive}
\newacronym{ISO}{ISO}{independednt system operator}

\newcommand{\0}{\bs 0}
\def\1{{\bs 1}}

\def\argmin{\mathop{\rm argmin}}
\newcommand{\col}{\mathrm{col}}

\def\diag{\mathop{\hbox{\rm diag}}}

\newcommand{\proj}{\mathrm{proj}}

\def\nc{\mathrm{N}}

\def\R{\mathbb{R}}
\def\N{\mathbb{N}}

\def\nbar{\skew2\bar n}

\def\xbar{{\bar x}}

\def\lbar{{\bar \l}}
\def\zbar{{\bar z}}

\def\I{i \in \mc{I}}
\def\diff#1{\textstyle \frac{\partial}{\partial #1}}
\def\x{\bs{x}}
\def\y{\bs{y}}
\def\z{\bs{z}}
\def\rhobs{\bs{\rho}}
\def\varsigmabs{\bs{\varsigma}}
\def\psibs{{\bs{\psi}}}
\def\w{\bs{\omega}}
\def\wbar{\bar{\bs{\omega}}}
\def\l{\bs{\lambda}}
\def\lbar{\bar{\bs{\lambda}}}
\def\z{\bs{z}}
\def\xbar{\bs{\bar{x}}}
\def\lbar{\bs{\bar{\lambda}}}
\def\zbar{\bs{\bar{z}}}
\def\L#1{\bs{L}_{#1}}
\def\I{i\in\mc{I}}
\def\K{k\in\mc{K}_i}
\def\nbar{{\bar{n}}}
\def\grad#1{ \nabla_{\!\! #1}}
\def\sigmabs{\bs{\sigma}}
\def\B{\mc{B}}

\DeclareSymbolFont{myletters}{OML}{ztmcm}{m}{it}
\DeclareMathSymbol{\uplambda}{\mathord}{myletters}{"15}

\def\QEDhereeqn{\eqno\let\eqno\relax\let\leqno\relax\let\veqno\relax\hbox{\QED}}
\def\QEDopenhereeqn{\eqno\let\eqno\relax\let\leqno\relax\let\veqno\relax\hbox{\QEDopen}}

\newtheorem{lem1}{Lemma}
\newtheorem{asm}{Assumption}
\theoremstyle{definition}
\newtheorem{defn1}{Definition}
\newtheorem{rem1}{Remark}

\begin{document}

\begin{frontmatter}
	
	\title{Continuous-time fully distributed generalized Nash equilibrium seeking for multi-integrator agents\thanksref{footnoteinfo}} 
	
	\thanks[footnoteinfo]{ This work is supported by NWO under project OMEGA (613.001.702) and by the ERC under  project COSMOS (802348). The material in this paper was partially presented at the 2020 European Control Conference.}
	
	\author[aff1]{Mattia Bianchi}\ead{m.bianchi@tudelft.nl},    
	\author[aff1]{Sergio Grammatico}\ead{s.grammatico@tudelft.nl},               

	\address[aff1]{Delft Center for Systems and Control, Delft University of Technology, The Netherlands}
	
		

\begin{abstract}
We consider strongly monotone
games with convex separable coupling constraints, played by dynamical
 agents, in a partial-decision information scenario.
We start by designing continuous-time fully distributed feedback controllers, based on consensus and primal-dual gradient dynamics, to seek a generalized Nash equilibrium in networks of single-integrator agents.
Our first solution adopts a fixed gain, whose choice requires the knowledge of some global parameters of the game. To relax this requirement, we  conceive a controller that can be tuned in a completely decentralized fashion,
thanks to the use of  uncoordinated integral adaptive weights. We further introduce algorithms specifically devised for generalized aggregative games. Finally, we adapt all our control schemes to deal with heterogeneous multi-integrator agents and, in turn, with nonlinear feedback-linearizable dynamical systems.   
For all the proposed dynamics, we show convergence to a variational equilibrium, by leveraging monotonicity properties and stability theory for projected dynamical systems. 
\end{abstract}

\end{frontmatter}


\section{Introduction}\label{sec:introduction}

Generalized games arise in several engineering applications, including demand-side management in the smart grid \citep{Saad2012}, charging scheduling of electric vehicles \citep{Grammatico2017} and communication networks \citep{Palomar_Eldar_facchinei_pang_2009}. 
These scenarios involve multiple autonomous decision makers, or agents; each agent aims at minimizing its individual cost function -- which depends on  its own action as well as on the actions of other agents -- subject to shared constraints.
Specifically, in many distributed control problems, the action of an agent  consists of the output of a dynamical system. For instance, in coverage maximization \citep{DurrStankovicJohanssonACC2011} and connectivity problems \citep{Stankovic_Johansson_Stipanovic_2012},  the agents are vehicles with some inherent dynamics, designed to optimize  inter-dependent objectives related to their positions; in electricity markets, the actions are represented by the power produced by some generators \citep{DePersisMonshizadeh_2019}; in optical networks, the costs are a function of the output powers of some  dynamical channels \citep{RomanoPavel2019}. 

In this context, the goal is to drive the physical processes
 to a desirable steady state, usually identified with a \gls{GNE},
 using only the  local information available to each agent. 
 One possibility is to exploit time-scale separation between the computation of a \gls{GNE} and setpoint tracking; yet, this solution is  typically economically inefficient and not robust 
 \citep{ZhangPapachristodoulouLiCDC2015}. 
 Alternatively, part of the recent literature focuses on the design of   distributed feedback controllers, to automatically steer a dynamical network to some (not known a priori) convenient operating point, while also ensuring closed-loop stability \citep{DePersisMonshizadeh_2019,DallAneseACC2015}. This paper fits in the latter framework. 

 In particular, we investigate \gls{GNE} seeking for multi-integrator agents, motivated by 
  robotics and mobile sensors applications \citep{FrihaufKrsticBasar_2012,Stankovic_Johansson_Stipanovic_2012}, where multi-integrator dynamics are commonly used to model elementary vehicles. 
  The study of this class of systems allows us to address \gls{GNE} problems for a variety of dynamical agents, linear or nonlinear, via feedback linearization (e.g.,  Euler--Lagrangian systems as in \cite{Deng_Liang_EL_NE_2019}).

 \emph{Literature review:} 
A variety of  algorithms has been proposed to seek a \gls{GNE} in a distributed way \citep{YiPavel2019,Yu_VanderSchaar_Sayed_2017,BorgensKanzow2018}  with a focus on aggregative games \citep{BelgioiosoGrammatico_L_CSS_2017,PaccagnanCDC2016,DePersisGrammaticoTAC2020}.
These works refer to (aggregative) games played in a full-decision information setting, where each agent can access the action of all the competitors (aggregate value), for example in  presence of a central coordinator that broadcasts the data to the network. 
Nevertheless, this is impractical in many applications, where the agents  only rely on local information. 
%

Instead, in this paper, we consider the so-called \emph{partial-decision information} scenario, where each agent holds an analytic expression for its cost  but is unable to evaluate the actual value, since it cannot access the strategies of all the competitors. To remedy the lack of knowledge, the agents agree on sharing some information with some  trusted neighbors over a communication graph. Based on  the data exchanged, each agent can estimate and asymptotically reconstruct the actions of all the other agents. This setup has been investigated for games without coupling constraints, resorting to gradient and consensus dynamics, both in discrete-time \citep{TatarenkoShiNedic_CDC2018, Koshal_Nedic_Shanbag_2016}, and continuous-time \citep{GadjovPavel2018,YeHu2017}. 
 Fewer works deal with generalized games \citep{Pavel2018,DengNian2019,PariseGentileLygeros_TCNS2020}.
 Moreover, all the results mentioned above consider static or single-integrator agents only. 
 Distributively driving  a network of more complex 
 physical systems to a \gls{NE} is still a relatively unexplored problem. 
With regard to aggregative games, a proportional integral feedback algorithm 
was developed in \cite{DePersisMonshizadeh_2019}  to seek a \gls{NE} in  networks of passive  second-order systems;  in \cite{Deng_Liang_EL_NE_2019} and \cite{Zhang_et_al_NL_NE_2019}, continuous-time gradient-based controllers were introduced for some classes of nonlinear dynamic. \cite{Stankovic_Johansson_Stipanovic_2012} addressed generally coupled cost games played by linear agents, via an extremum seeking approach;  \gls{NE} problems in systems of multi-integrator agents were studied by 
 \cite{RomanoPavel2019}.
Yet, none of  these works considers generalized games. 
 Despite the scarcity of results, the presence of shared constraints is a significant extension, which arises naturally  when the agents compete for common resources \cite[§2]{FacchineiKanzow2010}. 
However, dealing with coupling constraints in a distributed fashion is extremely challenging. All the results available resort to primal-dual reformulations \citep{Pavel2018,DengNian2019}, where the main technical complications are the loss of monotonicity properties of the original problem and the non-uniqueness of dual solutions. 

 \emph{Contributions:}
 Motivated by the above, we develop fully distributed continuous-time  controllers to seek a  \gls{GNE}  
  in networks of  multi-integrator agents.  
We focus on games 
with separable  coupling constraints, 
played under  partial-decision information. Our novel contributions are summarized as follows:
 \begin{itemize}[leftmargin=*]
 	\item {\it  Nonlinear coupling constraints:}
 	 	We introduce primal-dual projected-gradient controllers to drive
        single-integrator agents to a \gls{GNE}, with convergence guarantees under strong monotonicity and Lipschitz continuity of the game mapping. 
        In contrast with the existing fully distributed methods, we allow for arbitrary convex separable (not necessarily affine) coupling constraints. 
        Besides, our schemes are the only continuous-time fully distributed algorithms for generalized games (except for that in \cite{DengNian2019},  for aggregative games and specific  equality constraints only) (§3-4);
 	\item  { \it Adaptive \gls{GNE} seeking:}
 	We conceive the  first \gls{GNE} seeking algorithm that can be tuned in a fully decentralized way and without requiring any global information. Specifically, we extend the result in \cite{DePersisGrammatico2018} to generalized games and prove that convergence to an equilibrium can be ensured by adopting integral weights in place of a fixed, global, high-enough gain, whose choice would require the knowledge of the algebraic connectivity of the communication graph and of the Lipschitz and strong monotonicity constants of the game mapping
%
%
%
(§3-4);
 	\item {\it Generalized aggregative games:}
 	We propose controllers  for aggregative games with affine aggregation function, where the agents keep and exchange an estimate of the aggregation value only, thus reducing communication and computation cost.
 	Differently from the existing results, e.g., \cite{DengNian2019}, we can handle generic coupling constraints, thanks to a new variant of continuous-time dynamic tracking.  Furthermore, we develop an adaptive algorithm that requires no a priori  information and virtually no tuning (§5);
 	\item \textbf{\it Heterogeneous multi-integrator agents:}
 	We show how all  our controllers can be adapted to  solve \gls{GNE} problems  where \emph{each} agent is described by mixed-order  integrator dynamics, a class never considered before. Importantly, this allows us to address  games played by arbitrary   nonlinear agents with maximal relative degree, via feedback linearization. To the best of our knowledge, we are the first to study \emph{generalized} games with higher-order dynamical agents
(§6).
 \end{itemize}
To improve readability, the proofs are in the appendix.
Some preliminary results have been presented in \cite{Bianchi_ECC20_ctGNE}; the  novel contributions of this paper are:
 we consider adaptive controllers that can be tuned without need for any global information; we address a wider class of generalized games with nonlinear coupling constraints; we present algorithms for aggregative games, scalable with respect to the number of agents; we address the case of  mixed-order multi-integrators (instead of double-integrators); we provide a more extensive numerical analysis, including applications to networks of heterogeneous nonlinear systems. 
\emph{Basic notation:} See \cite{Bianchi_ECC20_ctGNE}.

\emph{Operator-theoretic definitions:}
An operator $\mc{F} : \R^n \rightarrow \R^n$ is
monotone ($\mu$-strongly monotone) if, for all $x, y \in \R^n$,  $(\mc{F}(x)-\mc{F}(y))^\top  (x-y) \geq 0$ $(\geq\mu \|x-y\|^2 )$.
For a  closed convex set $S\subseteq \R^n$,  $\proj_{S}:\R^n \rightarrow S$ is the Euclidean projection onto $S$;
$\nc_{S}: S \rightrightarrows \R^n:x\mapsto \{ v \in \R^n \mid \sup_{z \in S} \, v^\top (z-x) \leq 0  \}$ is the normal cone operator of $S$;
 $\mathrm{T}_S:S \rightrightarrows \R^n:x\mapsto \operatorname{cl}({\bigcup_{\delta>0}\frac{1}{\delta}(S-x)})$ is the tangent cone operator of $S$, where $\operatorname{cl}(\cdot)$ denotes the set closure.
The projection on the tangent cone of $S$ at $x$ is $\Pi_{S}(x,v):=\proj_{\mathrm{T}_S(x)}(v)=\lim_{\delta\rightarrow 0^+}\textstyle \frac{\proj_{S}(x+\delta v)-x}{\delta}$.
By Moreau's Decomposition Theorem \cite[Th.~6.30]{Bauschke2017}, $v=\proj_{\mathrm{T}_S(x)}(v)+\proj_{\nc_S(x)}(v)$ and 
$\proj_{\mathrm{T}_S(x)}(v)^\top \proj_{\nc_S(x)}(v)=0$, for any $v\in\R^n$.
%
%
%

\emph{Projected dynamical systems \citep{Cherukuri2016}:} 
Given an operator $\mc{F}:\R^n\rightarrow \R^n$ and a closed convex set $S\subseteq \R^n$, we consider the projected dynamical system
\begin{align}\label{eq:projdyn} \dot{x}=\Pi_S(x,\mc{F}(x)), \qquad x(0)=x_0\in S.
\end{align}
In \eqref{eq:projdyn}, the projection operator is possibly discontinuous on the boundary of 
$S$. If $\mc{F}$ is Lipschitz on $S$, the system \eqref{eq:projdyn} admits a unique global Carathéodory solution, i.e., there exists a unique absolutely continuous function $x:\R_{\geq 0}\rightarrow \R^n$ such that $x(0)=x_0$, $\dot{x}(t)=\Pi_S(x,g(x))$ for almost all $t$. Moreover, $x(t)\in S$ for all $t\geq 0$, as on the boundary of $S$ the projection operator 
restricts the flow of $\mc{F}$ such that the solution of \eqref{eq:projdyn} remains in $S$ (while $\Pi_S(x,\mc{F}(x))=\mc{F}(x)$ if $x\in \rm{int}(S)$).

\begin{lem1}\label{lem:minoratingproj}
	Let  $S\subseteq \R^q$ be a nonempty closed convex set. For any $y,y^\prime \in S$ and any $\xi \in \R^q$, it holds that $
	(y-y^\prime)^\top \Pi_{S}\left(y,\xi\right)\leq (y-y^\prime)^\top \xi$. 
	In particular, if $\ \Pi_S(y,\xi)=0$, then
	$
	(y-y^\prime)^\top\xi \geq 0
	$ (i.e., $\xi \in \nc_S(y)$).
	{\hfill $\square$}
\end{lem1}
\begin{pf*}{Proof.}
	By Moreau's theorem, $\left(\xi-\Pi_{S}(y,\xi)\right)\in \mathrm{N}_S(y)$; thus $\forall y,y^\prime\in S$,
	$
	(y^\prime-y)^\top(\xi- \Pi_{S}(y,\xi))\leq 0$. \hfill $\blacksquare$
\end{pf*}
\section{Mathematical Background}\label{sec:mathbackground}
We consider a group of agents $ \mc I:=\{ 1,\ldots,N \}$, where each agent $i\in \mc{I}$ shall choose its decision variable (i.e., strategy) $x_i$ from its local decision set $\textstyle \Omega_i \subseteq \R^{n_i}$. Let $x := \col(( x_i)_{i \in \mc I})  \in \Omega $ denote the stacked vector of all the agents' decisions, $\textstyle \Omega := \times_{i \in \mc I} \Omega_i \subseteq \R^n$ the overall action space and $\textstyle n:=\sum_{i=1}^N n_i$. 
The goal of each agent $i \in \mc I$ is to minimize its objective function $J_i(x_i,x_{-i})$, which depends  both on the local strategy $x_i$ and on the decision variables of the other agents $x_{-i}:= \col(( x_j)_{j\in \mc I\backslash \{ i \} })$.
Furthermore, we address \textit{generalized} games, 
where the coupling among the agents arises also via their feasible decision sets. In particular, we consider  separable coupling constraints, so that
 the overall feasible set is
$
\mc{X} :=\textstyle \Omega \cap\left\{x \in \R^{n} \mid g(x)\leq \0_m \right\},
$
where  $g:\R^n\rightarrow \R^m$, $g(x):=\textstyle \sum_{\I} g_i(x_i)$, and $g_i:\R^{n_i}\rightarrow \R^m$ is a private function of agent $i$.
%
The game is then  represented by $N$ inter-dependent optimization problems:
\begin{align} \label{eq:game}
\forall i \in \mc{I}:
\;  \underset{y_i \in \R^{n_i}}{\argmin}\; J_i(y_i,x_{-i})
\, \text{s.t. }   \;  (y_i,x_{-i}) \in \mc X.
\end{align}
The technical problem we consider in this paper is the computation of a \gls{GNE}, a joint action from which no agent has interest to unilaterally deviate.
\begin{defn1}
	\setlength{\belowdisplayskip}{0pt}
A  collective strategy $x^{*}=\operatorname{col}(\left(x_{i}^{*}\right)_{i \in \mathcal{I}})$ is a  generalized Nash equilibrium if, for all $i\in \mc{I}$,
\[ {x_{i}^{*} \in \underset{y_i}{\argmin} \, J_{i}\left(y_{i}, x_{-i}^{*}\right) \textup{ s.t. } (y_{i},x^*_{-i})\in \mc{X}.}    \QEDopenhereeqn
\]
\end{defn1}
Next, we formulate standard convexity and regularity
assumptions for the constraints and cost functions
(\citealp[Asm.~1]{KulkarniShanbag2012}; \citealp[Asm.~1]{Pavel2018}).
\begin{asm}\label{Ass:Convexity}
For each $i\in \mathcal{I}$, the set $\Omega_i$ is closed and convex;  $g_i$ is componentwise convex and twice continuously differentiable; $\mc{X}$
satisfies Slater's constraint qualification;  $J_{i}$ is continuously differentiable and the function $J_{i}\left(\cdot, x_{-i}\right)$ is convex for every $x_{-i}$.
{\hfill $\square$} \end{asm}
%

Under Assumption~\ref{Ass:Convexity},  $x^*$ is a \gls{GNE} of the game in \eqref{eq:game} if and only if there exist dual variables $\lambda_i^*\in \R^m $ such that the following  \vphantom{\gls{KKT}} \gls{KKT} conditions are satisfied, for all $\I$ \cite[Th.~4.6]{FacchineiKanzow2010}:
	\begin{align}\label{eq:KKT_i}
	\begin{aligned}
	{\0_{n_i}} & \in \nabla _{\! \! x_i} J_i(x_i^*,x_{-i}^*)+\diff{x_i}g_i(x_i^*)^\top \lambda_i^{*}+\nc_{\Omega_i}\left(x_i^{*}\right) 
		\\
	{\0_{m}} & \in-g(x^*)+\mathrm{N}_{\R_{\geq 0}}^{m}\left(\lambda_i^{*}\right).
	\end{aligned}
	\end{align}
Specifically, we focus on the subclass of \glspl{v-GNE} \cite[Def.~3.11]{FacchineiKanzow2010}, namely  \glspl{GNE} with equal dual variables, i.e.  $\lambda_i^*=\lambda^*\in\R^m$ for all $\I$, for which the \gls{KKT} conditions read as
\begin{subequations}\label{eq:KKT}
\begin{align} \label{eq:KKTa}
{\0_{n}} & \in F\left(x^{*}\right) +\diff{x}g(x^*)^\top \lambda^*+\mathrm{N}_{\Omega}(x^{*})
\\
\label{eq:KKTb}
{\0_{m}} & \in-g(x^*)+\mathrm{N}_{\R_{\geq 0}}^{m}\left(\lambda^{*}\right).
\end{align}
\end{subequations}
where $F$ is the \emph{pseudo-gradient} mapping of the game:
\begin{align}
\label{eq:pseudo-gradient}
F(x):=\operatorname{col}\left( (\grad{x_i} J_i(x_i,x_{-i})\right)_{i\in\mathcal{I}}) .
\end{align}
Variational equilibria enjoys important structural properties, such as economic fairness
\citep{FacchineiKanzow2010}.
For example, in electricity markets, the dual variables correspond to unitary prices charged for the use of the infrastructure  by an administrator  that aim at maximizing its revenue while ensuring certain operating conditions, and it is reasonable to assume that the administrator cannot charge discriminatory prices to different energy producers \citep{KulkarniShanbag2012}.
A sufficient condition for the existence and uniqueness of a \gls{v-GNE} is the strong monotonicity of the pseudo-gradient \cite[Th.~1, Rem.~1]{YiPavel2019}, which was always postulated in continuous-time \gls{NE} seeking under partial-decision information (\citealp[Asm.~2]{GadjovPavel2018}; \citealp[Asm.~3]{DengNian2019}). It implies strong convexity of the functions $J_i(\cdot, x_{-i} )$
 for any $x_{-i}$ \cite[Rem.~1]{TatarenkoShiNedic_CDC2018}, but not necessarily convexity of $J_i$ in the full argument.
\begin{asm}
	\label{Ass:StrMon}
The game mapping $F$  in \eqref{eq:pseudo-gradient} is:
\begin{itemize}[noitemsep, topsep=-4pt]
	\item[(i)]$\mu$-strongly monotone, for some $\mu>0$;
	\item[(ii)] $\theta_0$-Lipschitz continuous, for some $\theta_0>0$.\hfill $\square$
\end{itemize} 
%
 \end{asm}

%

\section{Fully distributed equilibrium seeking}\label{sec:distributedGNE}
In this section, we consider the game in  \eqref{eq:game}, where each agent is associated with  the following dynamical system:
\begin{equation}\label{eq:integrators}
\textstyle
\forall i \in \mathcal{I}: \quad \dot{x}_i=u_i, \ x_i(0)\in \Omega_{i}. 
\end{equation}
Our aim is to design the inputs $u_i \in \R^{n_i}$ to seek a \gls{v-GNE} in a fully distributed way. Specifically, each agent $i\in\mc{I}$ only knows its own feasible set $\Omega_i$, the portion $g_i$ of the coupling constraints,
and its own cost function $J_i$.
Moreover, the agents cannot access the strategies of all the competitors $x_{-i}$. 
Instead, each agent only relies on the information exchanged locally with some neighbors over a communication network $\mathcal G(\mc{I},\mc{E})$. The unordered pair $(i,j)$ belongs to the set of edges $\mc{E}$  if and only if agent $i$ and $j$ can exchange information. We denote  by $W\in \R^{N\times N}$ the symmetric adjacency matrix of $\mc{G}$, with $[W]_{i,j}>0$ if $(i,j)\in\mc{E}$, $[W]_{i,j}=0$ otherwise; $L$ the symmetric Laplacian matrix of $\mc{G}$; $\mc{N}_{i}:=\{j\mid(i,j)\in\mc{E}\}$  the set of neighbors of agent $i$. For ease of notation, we assume that the graph is unweighted, i.e., $[W]_{i,j}=1$ if $(i,j)\in \mc{E}$, but our results still hold for the weighted case.
\begin{asm}
	\label{Ass:Graph}
	The communication graph $\mathcal G (\mc{I},\mc{E}) $ is undirected and connected. 
{\hfill $\square$} \end{asm}
Our first algorithm is inspired by the discrete-time primal-dual gradient iteration in \cite[Alg.~1]{Pavel2018}.
To cope with the lack of knowledge, the general assumption for the partial-decision information scenario is that each agent keeps an estimate of all other agents’ actions \citep{Pavel2018,TatarenkoShiNedic_CDC2018}.
Let $\x^{i}:=\operatorname{col}((\x^{i}_j)_{j\in \mc{I}})\in \R^{n}$,  where $\x^{i}_{i}:=x_i$ and $\x^i_{j}$ is agent $i$'s estimate of agent $j$'s action, for all $j\neq i$; $\x^j_{-i}:=\col((\x_{j,\ell})_{\ell\in\mc{I}\backslash\{i\}})$.
Each agent also keeps
an estimate $\lambda_i\in \R_{\geq 0}^m$ of the dual variable and an auxiliary variable $z_i\in\R^m$ to allow for distributed consensus of the dual estimates. Our proposed dynamics are summarized in Algorithm~\ref{algo:1}, where $c>0$ is a \emph{global} fixed parameter (and $\theta$ is a constant  defined in Lemma~\ref{lem:LipschitzExtPseudo}). 
\begin{algorithm}[t] \caption{Constant gain} \label{algo:1}
    \textbf{Initialization}: set $c>\underline{c}:=\frac{(\theta_0+\theta)^{2}+4\mu\theta}{4\mu\uplambda_2(L)}$;
    $\forall \I$, set 	
	$  \bs{x}^i_{-i}(0)\in \R^{n-n_i}$,  $z_i(0)=\0_m$, $\lambda_{i}(0)\in \R^m_{\geq0}$;
		\vspace{0.5em}
		\\
	\textbf{Dynamics}: $\forall \I,$
	\vspace{-0.8em}
	\begin{align}
	\nonumber
	\begin{aligned}
	\dot{x}_i
	&=u_i=\Pi_{\Omega_{i}}(
	\begin{multlined}[t]
	x_i,
	-\nabla_{\!\! x_{i}} J_{i}(x_i,\bs{x}^{i}_{-i})-\diff{x_i}g_{i}(x_i)^{\top} \lambda_{i}\\- c{\textstyle \sum_{j \in \mathcal{N}_{i}}}
	(x_{i}-\bs{x}_{i}^{j}))
	\end{multlined}
	\\
	\dot{\bs{x}}^i_{-i}&=-c\textstyle\sum_{j \in \mathcal{N}_{i}} 
	(\bs{x}_{-i}^{i}-\bs{x}_{-i}^{j} ) 
	\\
	{\dot z}_i&=\textstyle \sum_{j \in \mathcal{N}_{i}} 
	(\lambda_{i}-\lambda_{j} )
	\\
	\dot\lambda_i&=
	\Pi_{\R_{\geq 0}^m} ( \lambda_i, g_i(x_i) -z_i-\textstyle \sum_{j\in\mc{N
		}_i}(\lambda_i-\lambda_j) )
	\end{aligned}
	\end{align}
	\vspace{-1.1em}
\end{algorithm}

 We note that the  agents exchange $\{ \bs{x}^i,\lambda_i \}$ with their neighbors only, therefore the controller can be implemented distributedly. Importantly,  each agent $i$ evaluates the partial gradient of its cost  $\grad{x_i} J_i$ on its local estimate $\x^i$, not on the actual joint strategy $x$.
 In steady state, the agents should agree on their estimates, i.e., $\bs{x}^i=\bs{x}^j, \ \lambda_i=\lambda_j$, for all $i,j \in \mc{I}$. This motivates the presence of consensual terms for both primal and dual variables. 
For any integer $q$, we denote  $\bs{E}_{q}:=\{\bs{y} \in \R^{N q}:\bs{y}=\1_{N} \otimes y, y\in \R^q\}$ the consensus subspace of dimension $q$, and $\bs{E}_{q}^\perp:=\{\boldsymbol{y} \in \mathbb{R}^{N q}:\left(\mathbf{1}_{N}^{\top} \otimes I_{q}\right) \boldsymbol{y}=\mathbf{0}_{q}\}$ its orthogonal complement;
Specifically, $\bs{E}_n$ and $\bs{E}_m$ are the action and multiplier consensus subspaces, respectively.
Moreover, $P_q:=\textstyle \frac{1}{N} \1_N \1_N^\top \otimes I_q$ is the projection matrix onto $\bs{E}_q$, i.e., $P_q \y=\proj_{\bs{E}_q}(\y)$,
 and 
 $P^{\perp}_q :=I_{Nq}-P_q$ the projection 
matrix onto the disagreement subspace $\bs{E}^{\perp}_q$.
\begin{algorithm}[t] \caption{Adaptive gains} \label{algo:2}
	\vspace{0.2em}
	\textbf{Initialization}: $\forall \I$, set $\gamma_i>0$, $\bs{x}^i_{-i}(0)\in \R^{n-n_i}$,  $k_i(0)\in \R$, $z_i(0)=\0_m$, $\lambda_{i}(0)\in \R^m$;
	\vspace{0.5em}
	\\
	\textbf{Dynamics}: $\forall \I$, 
	\vspace{-0.8em}
	\begin{align*} 
	\begin{aligned}
	\dot{x}_i& =u_i= \Pi_{\Omega_{i}}
	\begin{multlined}[t]
	(x_i,-\grad{x_i} J_{i}(x_{i}, \bs{x}_{-i}^{i})-\diff{x_i}g_{i}(x_i)^{\top} \lambda_{i}
	\\-  {\textstyle \sum_{j \in \mathcal{N}_{i}}} (k_i{{\rho}}^i_i-k_j{{\rho}}^i_i))
	\end{multlined}
	\\
	\dot{\bs{x}}^i_{-i}&=-\textstyle\sum _{j \in \mathcal{N}_{i}}
	(k_j{{\rho}}^i_{-i}-k_i{{\rho}}^j_{-i})
	\\
	\dot k_i&=\gamma_i \|{{\rho}}^i \|^2,  \qquad {{\rho}}^i=	\textstyle \sum_{ j \in \mathcal{N}_{i}} 
	(\bs{x}^{i}-\bs{x}^{j})
	\\
	{\dot z}_i&=\textstyle \sum_{j \in \mathcal{N}_{i}} \left(\lambda_{i}-\lambda_{j}\right) 
	\\
	\dot\lambda_i&=
	\Pi_{\R_{\geq 0}^m} (\lambda_i, g_i(x_i)- z_i -\textstyle \sum_{j\in\mc{N
		}_i}(\lambda_i-\lambda_j))
	\end{aligned}
	\end{align*}
	\vspace{-1.1em}
\end{algorithm}

While Algorithm~\ref{algo:1} is fully distributed, choosing the gain $c$  requires global knowledge about the graph $\mathcal G$, i.e., the algebraic connectivity, and about the game mapping, i.e., the strong monotonicity and Lipschitz constants. These parameters are  unlikely to be available locally in a network system.
To overcome this limitation and enhance  scalability, \cite{DePersisGrammatico2018}  proposed a controller where the communication gains are tuned online, thus relaxing the need for global information, for games without coupling constraints. Here we extend their solution to the \gls{GNE} problem.

Our proposed controller is given in Algorithm~\ref{algo:2}.
For all $i\in \mathcal I$, $k_i$ is the adaptive gain of agent $i$, 
$\gamma_i>0$ is an arbitrary local constant and  $\rho^i:=\col((\rho^i_j)_{j\in\mc{I}})$.
We emphasize that 
	Algorithm~\ref{algo:2}  allows for a fully uncoupled tuning: each agent chooses \emph{locally} the initial conditions and the  parameter $\gamma_i$, independently of the other agents and without  need for coordination or global knowledge. 
\begin{rem1}
 Algorithm~\ref{algo:2} uses second order information, as each agent sends the quantity $\rho^i$, which depends on the estimates of its neighbors. In  case of delayed communication, this means dealing with twice the transmission latency with respect to a controller that exploits first order information only, e.g., Algorithm~\ref{algo:1}. In a discrete-time setting,  a sampled version of Algorithm~\ref{algo:2} can be implemented by allowing  the agents to  communicate twice per iteration, a common assumption for \gls{GNE} seeking on networks \citep{Pavel2018,GadjovPavel_aggregative_2019}.
	\hfill $\square$
\end{rem1}
 To rewrite the closed-loop dynamics in Algorithms~\ref{algo:1}, \ref{algo:2}  in compact form, let us define $\x:=\col((\x^i)_{i\in\mc{I}})$ and, as in \cite[Eq.~11]{GadjovPavel2018}, for all $i \in \mc{I}$,
\begin{align}
\mathcal{R}_{i}:&=
[ \begin{array}{llll}{{\0}_{n_{i} \times n_{<i}}} & {I_{n_{i}}} & {\0_{n_{i} \times n_{>i}}}\end{array} ], 
\end{align}
where $n_{<i}:=\sum_{j<i,j \in \mathcal{I}} n_{j}$, $n_{>i}:=\sum_{j>i, j \in \mathcal{I}} n_{j}$;  let also $\mathcal{R}:=\operatorname{diag}\left((\mathcal{R}_{i})_{i \in \mathcal{I}}\right)$. In simple terms, $\mathcal R _i$ selects the $i$-th $n_i$ dimensional component from an $n$-dimensional vector, i.e., 
$\mathcal{R}_{i} \bs{x}^{i}=\bs{x}_{i}^{i}=x_i$ and 
$x=\mathcal{R} \bs{x}$.

Let  $\bs \lambda:=\operatorname{col}((\lambda_i)_{i\in\mathcal I})$,  $\bs z:=\operatorname{col}((z_i)_{i\in\mathcal I})$, 
$\bs{\Omega}:=\{\x\in\R^{nN}\mid \mc{R}\x\in \Omega\}$,
$\bs{g}(x):=\col((g_i(x_i)_{\I}))$, 
$\mathbb{G}(x):=\diff{x}\bs{g}(x)=\diag ((\diff{x_i}g_i(x_i))_{\I})$,
$\bs{k} :=\operatorname{col}((k_{i})_{\I})$,
${\bs{\rho}} :=\operatorname{col}(({{\rho}}^{i})_{\I})$,
$K :=\operatorname{diag}((k_{i}I_n)_{\I})$,
$D({\bs{\rho}}) :=\operatorname{diag} (({{\rho}}^{i})_{\I})$, $\Gamma :=\operatorname{diag}((\gamma_{i})_{\I})$, 
and, for any integer $q>0$, $\bs{L}_{q}:=L \otimes I_{q}$.  Furthermore, we define the \emph{extended pseudo-gradient} mapping  $\bs{F}$ as: 
\begin{align}
\label{eq:extended_pseudo-gradient}
\bs{F}(\bs{x}):=\operatorname{col}( (\grad{x_{i}} J_{i}(x_{i}, \bs{x}_{-i}^{i}))_{i \in \mathcal{I}}).
\end{align}
Therefore, Algorithm~\ref{algo:1}, in compact form, reads as 
\begin{subequations}
	\label{eq:dynamics:CG}
	\begin{align} 
	\dot { \bs{x}}&=
	\Pi_{\bs{\Omega}} \left(\x, -\mc{R}^\top (\bs{F} (\bs{x})+\mathbb{G} (\mc{R} \bs{x})  ^\top \bs{\lambda})-c \bs{L}_{n} \bs{x} \right)
	\label{eq:dynamics:CG:a}
	\\
	\label{eq:dynamics:CG:b}
	\dot{\bs{z}}&= 
	\bs{L}_{m} \bs{\lambda}
	\\
	\label{eq:dynamics:CG:c}
	\dot{\bs{\lambda}}&= \Pi_{\R_{\geq 0}^{N m}}(\bs{\lambda}, \bs{g}(\mc{R}\bs{x}) - \bs{z}-\bs{L}_m\bs{\lambda} ),
	\end{align}
\end{subequations}
and Algorithm~\ref{algo:2} as
\begin{subequations}
	\label{eq:dynamics:AG}
	\begin{align}
	\label{eq:dynamics:AG:1}
	\dot { \x}&= 
	\Pi_{\bs{\Omega}}   ( \x, -\mc{R}^{\! \top}\!(\bs{F} (\bs{x})+\mathbb{G} (\mc{R} \bs{x}) ^{\!\top} \bs{\lambda})   - \bs{L}_{n} K \bs{\rho} ) 
	\\
	\label{eq:dynamics:AG:2}
	\dot{\bs k}&= D\left( {\bs{\rho}} \right)^\top \left( \Gamma \otimes I_n\right){\bs{\rho}}, \qquad {\bs{\rho}}=\bs{L}_n \bs{x} 
	\\
	\label{eq:dynamics:AG:4}
	\dot{\bs{z}}&= 
	\bs{L}_{m} \bs{\lambda}
	\\
	\label{eq:dynamics:AG:5}
	\dot{\bs{\lambda}}&= \Pi_{\R_{\geq 0}^{N m}}(\bs{\lambda},\bs{g}(\mc{R}\bs{x}) - \bs{z} -\bs{L}_m\bs{\lambda}).
	\end{align}
\end{subequations}

\section{Convergence analysis}
In this section, we show the convergence of our dynamics to a \gls{v-GNE}. 
We focus on the analysis of Algorithm~\ref{algo:2}, which presents more technical difficulties; the convergence of Algorithm~\ref{algo:1} can be shown analogously.

We start by noting an invariance property of our controllers, namely that if $\z(0)\in\bs{E}^\perp_m$ (for instance, $\z(0)=\0_m$), then $\z\in\bs{E^\perp_m}$ along any solution of \eqref{eq:dynamics:AG}, by \eqref{eq:dynamics:AG:4}.
The next lemma relates a class of equilibria of the system in \eqref{eq:dynamics:AG} to the \gls{v-GNE} of the game in \eqref{eq:game}. 
\begin{lem1}\label{lem:equilibria-vGNE-adaptive}
	Under Assumptions~\ref{Ass:Convexity}, \ref{Ass:StrMon}, \ref{Ass:Graph},  the following statements hold:
	\begin{enumerate}[leftmargin=*]
		\item[i)] Any equilibrium point $\operatorname{col}\left( \bar{\bs{x}},\bar{\bs k},\bar{\bs z},\bar{\bs \lambda}\right) $ of \eqref{eq:dynamics:AG} with  $\bs{\bar{z}}\in\bs{E}_m^\perp$ is such that $\bar{\bs{x}}=\1_{N} \otimes x^{*}$, $\bar{\bs \lambda}=\1_{N} \otimes \lambda^{*}$, where  the pair ($x^*,\lambda^*$) satisfies the KKT conditions in \eqref{eq:KKT}, hence $x^*$ is the \gls{v-GNE} of the game in \eqref{eq:game}.
		\item[ii)] 
		The system \eqref{eq:dynamics:AG} admits at least one equilibrium $\operatorname{col}\left( \bar{\bs{x}},\bar{\bs k},\bar{\bs z},\bar{\bs \lambda}\right) $ with $\bs{\bar{z}}\in\bs{E}^\perp_m$.
\hfill$\square$ 
	\end{enumerate} 
\end{lem1}
We remark that in Algorithm~\ref{algo:2} (or \ref{algo:1}) each agent $i$ evaluates the quantity $\grad{x_i} J_i$ in its local estimate $\x^i$, not on $x$. As a consequence, the operator $\mc{R}^\top \bs{F}$ is very rarely monotone, 
 even under strong monotonicity of the game mapping $F$. The loss of monotonicity is indeed the main technical difficulty arising in the partial-decision information scenario.
Following  \cite{Pavel2018}, \cite{DePersisGrammatico2018}, we deal with this issue by leveraging a  restricted monotonicity property, which can be guaranteed for any game satisfying  Assumptions \ref{Ass:Convexity}-\ref{Ass:Graph}, without additional hypotheses, as shown in the next lemmas.
The proof relies on the decomposition of $\bs{x}$ along the consensus space  $\bs{E}_{n}$, where $\bs{F}$ is strongly monotone, and the disagreement space $\bs{E}_{n}^{\perp}$, where $\bs{L}_n$ is strongly monotone.
	
\begin{lem1}[{\citealp[Lem.~3]{Bianchi_ECC20_ctGNE}}]\label{lem:LipschitzExtPseudo}
	Let Assumption~\ref{Ass:StrMon} hold. Then, the  mapping $\bs{F}$ in \eqref{eq:extended_pseudo-gradient} is $\theta$-Lipschitz continuous, for some $\theta\in [\mu,\theta_0]$.
	{\hfill $\square$} \end{lem1}

\begin{lem1}\label{lem:strongmon_adaptive}
	Let Assumptions~\ref{Ass:StrMon}, \ref{Ass:Graph} hold, and let 
	\begin{align}
	\label{eq:M2}
	{M}_1:= \hspace{-0.3em}\begin{bmatrix}{\frac{\mu}{N}} & \ {-\frac{\theta_0+\theta}{2\sqrt{N}}} \\ {-\frac{\theta_0+\theta}{2\sqrt{N}}} & \ {k^*\uplambda_{2}(L)^2 -\theta}\end{bmatrix}\hspace{-0.3em}, 
	\; \underbar{k} :=\textstyle \frac{(\theta_0+\theta)^{2}+4\mu\theta}{4\mu\uplambda_2(L)^2}.
	\end{align}
	For any $k^*>\underbar{k}$ and $K^*=I_{Nn}k^*$, for any $\bs{x}\in R^{Nn}$ and any $\bs{y}\in \bs{E}_{n}$, it holds that ${M}_{2} \succ 0$ and also that
\[
	 \begin{aligned}[b]
	\left(\bs{x}-\bs{y}\right)^\top(\mathcal{R}^\top\left(\bs{F}(\bs{x})-\bs{F}\left(\bs{y}\right)\right)+ \bs{L}_n K^* \bs{L}_n \left(\bs{x}-\bs{y}\right))
	\\ \geq \uplambda_{\min }({M}_{1})\left\|\bs{x}-\bs{y}\right\|^{2}.
	\end{aligned}  
	\QEDopenhereeqn
\]
\end{lem1}
We can now present the main result of this section. 

\begin{thm}[Convergence of Algorithm~\ref{algo:2}]\label{th:main2}
	Let Assumptions \ref{Ass:Convexity}, \ref{Ass:StrMon}, \ref{Ass:Graph} hold.
	For any initial condition in $S=\bs{\Omega}\times\R^{N}\times \bs{E}^\perp_m\times \R^{mN}_{\geq 0}$, the system in \eqref{eq:dynamics:AG} has a unique Carathéodory solution, which belongs to $S$ for all $t\geq 0$.  The solution converges to an equilibrium $\operatorname{col}\left( \bar{\bs{x}},\bar{\bs k},\bar{\bs z},\bar{\bs \lambda}\right) $, with  $\bar{\bs{x}}=\1_{N} \otimes x^{*}$, $\bar{\bs \lambda}=\1_{N} \otimes \lambda^{*}$, and
    $(x^*,\lambda^*)$ satisfies the KKT conditions in \eqref{eq:KKT}, hence  $x^*$ is the  \gls{v-GNE} of the game in \eqref{eq:game}.
	{\hfill $\square$}
\end{thm}

A similar result holds also for the the dynamics in \eqref{eq:dynamics:CG}.
\begin{thm}[Convergence of Algorithm~\ref{algo:1}]\label{th:main1}
	Let Assumptions~\ref{Ass:Convexity}, \ref{Ass:StrMon}, \ref{Ass:Graph} hold.
	Let $c>\underline{c}$, with
$\underline{c}$ as in Algorithm~\ref{algo:1}.
	For any initial condition in $S=\bs{\Omega}\times \bs{E}^\perp_m\times \R^{mN}_{\geq 0}$ the system in \eqref{eq:dynamics:CG} has a unique Carathéodory solution, which belongs to $S$ for all $t\geq 0$.  The solution converges to an equilibrium $\operatorname{col}\left( \bar{\bs{x}},\bar{\bs z},\bar{\bs \lambda}\right) $, with  $\bar{\bs{x}}=\1_{N} \otimes x^{*}$, $\bar{\bs \lambda}=\1_{N} \otimes \lambda^{*}$, and
	  $(x^*,\lambda^*)$ satisfies the KKT conditions in \eqref{eq:KKT}, hence $x^*$ is the  \gls{v-GNE} of the game in \eqref{eq:game}.
	{\hfill $\square$}
\end{thm}

\begin{rem1}\label{rem:dualizedlocalconstraints}
As for Euclidean projections, evaluating  $\Pi_{\Omega_{i}}(x,v)$ can be computationally expensive. If,  for some $\I$ and some twice continuously differentiable mapping $g_i^{\textnormal{loc}}$, $\Omega_i=\{x_i\in \R^{n_i} \mid g_i^{\textnormal{loc}}(x_i)\leq \0_p \}$,  then the following alternative updates can be used in Algorithm~\ref{algo:1} (and similarly in	 Algorithm~\ref{algo:2}):
\[
\begin{aligned}
\dot{x}_i& =
\begin{multlined}[t]
-\nabla_{\!\! x_{i}} J_{i}(x_i,\bs{x}^{i}_{-i})-\diff{x_i}g_{i}(x_i)^{\top} \lambda_{i} \\ \quad -\diff{x_i}g_{i}^{\textnormal{loc}}(x_i)^{\top} \lambda^{\textnormal{loc}}_{i} - c{\textstyle \sum_{j \in \mathcal{N}_{i}}}	(x_{i}-\bs{x}_{i}^{j})
\end{multlined} 
\\
\dot{\lambda}_i^{\textnormal{loc}} &=\Pi_{\R_{\geq 0}^p} ( \lambda_i, g_i^{\textnormal{loc}}(x_i)).
\end{aligned}
\]
In simple terms, the local constraints are dualized like the coupling constraints; but the corresponding dual variables are managed locally. The drawback of this primal-dual approach is that the satisfaction of the local constraints can only be ensured asymptotically.
\hfill $\square$
\end{rem1}

\section{Generalized aggregative games}
 \label{sec:aggregative}
\begin{algorithm}[t] \caption{Constant gain (aggregative games)} \label{algo:3}
	\textbf{Initialization}: set $c> \underline{c}:= \textstyle \frac{(\tilde{\theta}_\sigma)^{2}}{4\mu\uplambda_2(L)}$;
	$\forall \I$, set 	
	$\varsigma_i=\0_{\bar{n}}$, $z_i(0)=\0_m$, $\lambda_{i}(0)\in \R^m_{\geq0}$;
	\vspace{0.5em}
	\\
	\textbf{Dynamics}: $\forall \I,$
	\vspace{-0.8em}
	\begin{align}
	\nonumber
	\begin{aligned}
	\dot{x}_i
	&=u_i=\Pi_{\Omega_{i}}(
	\begin{multlined}[t]
	x_i, -\grad{x_{i}} f_{i}(x_{i}, {\sigma}^i) -\diff{x_i}g_{i}(x_i)^{\top} \lambda_{i}
	\\
	-c{B_i^\top}\textstyle\sum _{j \in \mathcal{N}_{i}}
	({\sigma}^i-{\sigma}^j))
	\end{multlined}
	\\
	\dot{\varsigma}_i&=-c\textstyle\sum_{j \in \mathcal{N}_{i}}  ({\sigma}^i-{\sigma}^j), \qquad {\sigma}^i={\psi_i(x_i)}+\varsigma_i
	\\
	{\dot z}_i&=\textstyle \sum_{j \in \mathcal{N}_{i}} (\lambda_{i}-\lambda_{j} )
	\\
	\dot\lambda_i&=
	\Pi_{\R_{\geq 0}^m} (  \lambda_i, {g_i(x_i) -z_i}-\textstyle \sum_{j\in\mc{N}_i}(\lambda_i-\lambda_j) )
	\end{aligned}
	\end{align}
	\vspace{-1.1em}
\end{algorithm}
 \FloatBarrier
\begin{algorithm}[t] \caption{Adaptive gains (aggregative games)} \label{algo:4}
	\vspace{0.2em}
	\textbf{Initialization}: $\forall \I$, set $\gamma_i>0$, $\varsigma_i=\0_{\bar{n}}$, $k_i(0)\in\R$, {$z_i(0)=\0_m$}, $\lambda_{i}(0)\in \R^m_{\geq0}$;  
	\vspace{0.5em}
	\\
	\textbf{Dynamics}: $\forall \I$,
	\vspace{-0.8em}
	\begin{align}
	\nonumber
	\begin{aligned}
	\dot{x}_i&=u_i=	\Pi_{\Omega_{i}}(
	\begin{multlined}[t]
	x_i, - \grad{x_{i}} f_{i}(x_{i}, {\sigma}^i) {-\diff{x_i}g_{i}(x_i)^{\top} \lambda_{i}}
	\\
	{-B_i^\top}\textstyle\sum _{j \in \mathcal{N}_{i}}(k_i\rho^i -k_j\rho^j) )
	\end{multlined}
	\\
	\dot{\varsigma}_i&=-\textstyle\sum_{j \in \mathcal{N}_{i}}  (k_i{\rho}^i-k_j{\rho}^j)  \qquad \qquad	{\sigma}^i={\psi_i(x_i)}+\varsigma_i
	\\
	\dot k_i&=\gamma_i \|{\rho}^i \|^2  \qquad \qquad \qquad \quad {\rho}^i=	\textstyle \sum_{j \in \mathcal{N}_{i}}  \left({\sigma}^{i}-{\sigma}^{j}\right) 
	\\
	{\dot z}_i&=\textstyle \sum_{j \in \mathcal{N}_{i}}(\lambda_{i}-\lambda_{j} )
	\\
	\dot\lambda_i&=
	\Pi_{\R_{\geq 0}^m} (  \lambda_i, {g_i(x_i) -z_i}-\textstyle \sum_{j\in\mc{N}_i}(\lambda_i-\lambda_j) )
	\end{aligned}
	\end{align}
	\vspace{-1.1em}
\end{algorithm}
In this section, we study  aggregative games, where the cost function of each agent depends only on the local action  and on an aggregation value  $\psi(x):=\frac{1}{N}\textstyle{ \sum_{i\in\mc{I}} } \psi_i(x_i)$, where $\psi_i:\R^{n_i}\rightarrow \R^\nbar$, for all $\I$. 
It follows that, for each $\I$, there is a function $f_i:\R^{n_i}\times \R^{\bar{n}}\rightarrow {\R}$ such that the original cost function $J_i$ in \eqref{eq:game} can be written as
\begin{align}\label{eq:Jf_aggregative}
J_i(x_i,x_{-i})=f_i(x_i,\psi(x)).
\end{align}
In particular, we focus on games with affine aggregation functions, where, for all $\I$, 
$\psi_i(x_i)=B_ix_i+d_i$, for some $B_i\in\R^{\nbar \times n_i}$, $d_i\in\R^{\nbar}$.
As a special case, this  class  includes the  common (weighted) average aggregative games.

Since an aggregative game is only a particular instance of the game in \eqref{eq:game},
Algorithms~\ref{algo:1}-\ref{algo:2} could still be used to drive the system \eqref{eq:integrators} to the \gls{v-GNE}. This would require each agent $i$ to keep (and exchange) an estimate of all other agents' actions, i.e., a vector of $n-n_i$ components; however, the cost of each agent is only a function of the aggregation value $\operatorname{\psi}(x)$, whose dimension $\bar{n}$ is independent of the number of agents. To reduce the communication and computation burden, we introduce two distributed controllers that are scalable with the number of agents, specifically designed to seek a \gls{v-GNE} in aggregative games. 
Our proposed dynamics are shown in Algorithms~\ref{algo:3} and \ref{algo:4}.

Since the agents rely on local information only, they do not have access to  the actual value of the aggregation $\psi(x)$. Hence, we embed each agent with an auxiliary error variable $\varsigma_i\in \R^{\bar{n}}$, which is an estimate of $\psi(x)-\psi_i(x_i)$. Each agent aims at asymptotically reconstructing the true aggregation value, based on the information received from its neighbors. We use the notation
\[
\grad{x_{i}} f_{i}(x_{i}, {\sigma}^i):=
\begin{multlined}[t]
\grad{y} f_{i}(y, {\sigma}^i)|_{y=x_i} + \textstyle \frac{1}{N}
B_i^\top\grad{y}f_i(x_i,y)|_{y={\sigma}^i}.
\end{multlined}
\]
We note that, in Algorithms~\ref{algo:3} and \ref{algo:4}, the agents exchange the quantities ${\sigma}^i\in \R^{\bar{n}}$, instead of  the variables $\bs{x}^i, \rho^i  \in\R^{n}$, like in Algorithms \ref{algo:1} and \ref{algo:2}. 
Let  $\bs{\sigma}:=\operatorname{col}(({\sigma}^i)_{i\in\mathcal I})$. We define the \emph{extended pseudo-gradient} mapping  $\bs{\tilde{F}}$ as 
\begin{align}
\label{eq:extended_pseudo-gradient_agg}
\bs{\tilde{F}}(x,\bs{\sigma}):=\operatorname{col}\big(\left(\grad{x_{i}} f_i(x_{i}, {\sigma}^{i})\right)_{i \in \mathcal{I}} \big).
\end{align}
\begin{asm}\label{Ass:LipschitzExtPseudoagg}
	The mapping $\bs{{\tilde F}}$ in \eqref{eq:extended_pseudo-gradient_agg} is $\tilde{\theta}$-Lipschitz continuous, for some $\tilde{\theta}>0$.
		Hence, $\bs{\tilde{F}}(x,\cdot)$ is $\tilde{\theta}_\sigma$-Lipschitz continuous, for some $\tilde{\theta}_\sigma \in(0,\tilde{\theta}]$,  $\forall x\in\R^{{n}}$. \hfill $\square$
\end{asm}
%
{Assumption~\ref{Ass:LipschitzExtPseudoagg} is standard  (\citealp[Asm.~4]{GadjovPavel_aggregative_2019}; \citealp[Asm.~3]{Koshal_Nedic_Shanbag_2016}) and can be shown to  hold under  Assumption~\ref{Ass:StrMon} if the matrix $[B_1 \, \dots \, B_N]$ is full row rank, e.g., for average aggregative games.
%

By defining $\bs{\varsigma}:=\col((\varsigma_i)_{\I})$, $K:=\diag((k_i I_\nbar)_{\I})$, {$\bs{\psi}(x):=\col((\psi_i(x_i))_{\I})$, $B:=\diag((B_i)_{\I})$},   the dynamics in Algorithms  \ref{algo:3} and \ref{algo:4} read, in compact form, as
\begin{subequations}
	\label{eq:dynamics:CGagg}
	\begin{align}\label{eq:dynamics:CGagg:a}
	 \dot { {x}}&= \Pi_{\Omega} (x, -\bs{\tilde{F}} (x,\bs{\sigma})-{\mathbb{G} (x) ^{\!\top} \bs{\lambda}} -c{B^{\!\top}}\bs{L}_{\bar{n}} \bs{\sigma})
	\\
	\label{eq:dynamics:CGagg:b}
	\dot{\bs\varsigma}&=-c\bs{L}_{\bar{n}} \bs{\sigma}, \qquad 	\bs{\sigma}={\bs{\psi}(x)}+\bs\varsigma   
	\\
	\label{eq:dynamics:CGagg:c}
	\dot{\bs{z}}&= 
	\bs{L}_{m} \bs{\lambda}
	\\
	\label{eq:dynamics:CGagg:d}
	\dot{\bs{\lambda}}&= \Pi_{\R_{\geq 0}^{N m}}\left( \l, {\bs{g}(x) -\bs{z} }-\bs{L}_{m} \bs{\lambda}\right),
	\end{align}
\end{subequations}
and
 \begin{subequations}
 	\label{eq:dynamics:AGagg}
 	\begin{align}\label{eq:dynamics:AGagg:a}
 	\dot { {x}}&= \Pi_{\Omega} ( x, -\bs{\tilde{F}} (x,{\bs\sigma})\!-{\mathbb{G} (x) ^{\!\top} \bs{\lambda}} -{B^{\!\top}} \bs{L}_{\nbar}K\rhobs)
 	\\
 	\label{eq:dynamics:AGagg:b}
 	\dot{\bs\varsigma}&=-\L{\nbar}K \rhobs, \qquad    \;	{\bs\sigma}={\bs{\psi}(x)}+\bs{\varsigma }\qquad
 	\\
 	\label{eq:dynamics:AGagg:d}
 	\dot{\bs k}&= D\left( {\bs{\rho}} \right)^\top \left( \Gamma \otimes I_{\bar{n}}\right){\bs{\rho}}, \qquad {\bs{\rho}}= \L{\nbar}{\bs\sigma} \quad 
 	\\
 	\label{eq:dynamics:AGagg:e}
 	\dot{\bs{z}}&= 
 	\bs{L}_{m} \bs{\lambda}
 		\\
 	\label{eq:dynamics:AGagg:f}
 	\dot{\bs{\lambda}}&=\Pi_{\R_{\geq 0}^{N m}}\left(\l, {\bs{g}(x) -\bs{z}}-\bs{L}_{m} \bs{\lambda} \right) ,
 	\end{align}
 \end{subequations}
 respectively.
 We note that
only if  the estimates of all the agents coincide with the actual value, i.e., $\bs{\sigma}=\1_N\otimes \psi(x)$, we can conclude that $\bs{\tilde{F}}(x,\bs{\sigma})=F(x)$, $F$ as in \eqref{eq:pseudo-gradient}.
 
 \begin{rem1}\label{rem:tracking}
 {From the updates in \eqref{eq:dynamics:CGagg:b} (or \eqref{eq:dynamics:AGagg:b}), we can infer an invariance property of the closed-loop system \eqref{eq:dynamics:CGagg} (or \eqref{eq:dynamics:AGagg}), namely that, at any time, $\textstyle \frac{1}{N}\sum_{\I}  \varsigma_i=\0_\nbar$,  and thus  $\textstyle \frac{1}{N}\sum_{\I}  \sigma_i=\psi(x)$ (or equivalently, $P_\nbar \sigmabs =\1_N\otimes \psi(x)$), provided that 
$\varsigmabs(0)=\0_{N\nbar}$. In fact, the dynamics of $\sigma_i$ in Algorithm~\ref{algo:3} can be regarded as a continuous-time dynamic tracking for the time-varying quantity $\psi(x)$, i.e., $\sigma_i(0)=\psi_i(x_i(0))$ and
 	\begin{align}\label{eq:tracking}
 	\dot{\sigma}_i= -c\textstyle\sum_{j \in \mathcal{N}_{i}}  ({\sigma}^i-{\sigma}^j) + \frac{d}{dt}(\psi_i(x_i)). \end{align}
We emphasize that in Algorithm~\ref{algo:3} there is no agent that knows the quantity $\psi(x)$. This is the main difference with respect to Algorithm~\ref{algo:1}, where the consensus of the estimates works instead as a leader-follower protocol. If the actions $x$ are constant, the dynamics in \eqref{eq:tracking} reduce to a standard average consensus algorithm and ensure that $\sigmabs\rightarrow \1_N\otimes \psi(x)$ exponentially, under  Assumption~\ref{Ass:Graph}. 
Therefore, when the action dynamics \eqref{eq:dynamics:CGagg:a} are \gls{ISS} with respect to the estimation error, convergence can be ensured via small-gain arguments (for $c$ big enough) -- a similar approach was used in \cite{DengNian2019}. However, in the presence of generic coupling constraints (even affine), this robustness cannot be guaranteed; to still ensure convergence, we  design an extra consensual term for the action updates, i.e. $c B^\top \L{\nbar}\sigmabs$.
Furthermore, via the \emph{error} variable $\varsigmabs$, we avoid studying the discontinuous dynamics in \eqref{eq:tracking}. 
 We finally note that we consider games with affine $\psi$ (a broader class than \cite{GadjovPavel_aggregative_2019}), but nonlinear aggregation functions are also  studied \citep{DengNian2019,Deng_Liang_EL_NE_2019,Zhang_et_al_NL_NE_2019}. However,  \cite{Deng_Liang_EL_NE_2019} and \cite{Zhang_et_al_NL_NE_2019} postulate strong monotonicity of an augmented operator, a condition much more restrictive than our  Assumption~\ref{Ass:StrMon}(i)  \cite[Rem. 2]{Deng_Liang_EL_NE_2019}; instead, the approach in \cite{DengNian2019} is not suitable to deal with generic coupling constraints, as discussed above. 
 }
\hfill $\square$
\end{rem1}
 By leveraging the invariance property  in Remark~\ref{rem:tracking}, we can obtain a refinement of Lemma~\ref{lem:strongmon_adaptive}.


\begin{lem1}\label{lem:strongmon_adaptive_agg}
	Let Assumptions~\ref{Ass:StrMon}(i), \ref{Ass:Graph}, \ref{Ass:LipschitzExtPseudoagg} hold, and let 
	\begin{align}\label{eq:M4}
	{M}_{2}=
	\begin{bmatrix}
	{\mu} & \ {-\frac{\tilde{\theta}_\sigma}{2}}
	\\
	{-\frac{\tilde{\theta}_\sigma}{2}} & \ {k^* \lambda_{2}(L)^2}
	\end{bmatrix}, \quad  
	\underline{k}=\textstyle \frac{\tilde{\theta}_\sigma^{2}}{4\mu  \uplambda_{2}(L)^2}
	\end{align}
	  For any $k^*>\underline{k}$ and $K^*=I_{N\nbar}k^*$, for any $(x,\bs{\sigma})$ such that $P_{\nbar}\bs{\sigma}=P_{\nbar} \psibs (x)$
	and any $(x^\prime,\bs{{\sigma}}^\prime)$ such that ${\bs{\sigma}}^\prime=P_{\nbar} \psibs(x^\prime)=\1_N\otimes \psi(x^\prime)$, it holds that ${M}_{2} \succ 0$,  and   that
   \begin{flalign*}
   &&
   \begin{aligned}[b]
   &{(x-x^\prime)^\top(\bs{\tilde{F}}(x,{\bs\sigma})-	\bs{\tilde{F}}(x^\prime,{{\bs\sigma}^\prime})) }
   \\
	+&{ ({\bs\sigma}-{{\bs\sigma}^\prime})^\top\L{\nbar}K^*\L{\nbar}({\bs\sigma}-{{\bs\sigma}^\prime})} 
   	\\
   	&
	\geq 
	\uplambda_{\min }({M}_{2})
	\left\|
	\col \left(
	x-x^\prime,
	{\bs\sigma}
	-{\1_N\otimes \psi(x) }\right)
	\right\|^2.
	\end{aligned} && \square
	\end{flalign*}
\end{lem1}
Next, we exploit Lemma~\ref{lem:strongmon_adaptive_agg} to  prove the  convergence of
Algorithm~\ref{algo:4}. An analogous result holds for Algorithm~\ref{algo:3}.
\begin{thm}[Convergence of Algorithm~\ref{algo:4}]\label{th:main4}
	Let Assumptions~\ref{Ass:Convexity}, \ref{Ass:StrMon}(i), \ref{Ass:Graph}, \ref{Ass:LipschitzExtPseudoagg} hold.
	Then, for any initial condition in {$S={\Omega}\times \bs{E}^\perp_\nbar \times\R^N\times \bs{E}^\perp_{m} \times \R^{mN}_{\geq 0}$} the system in \eqref{eq:dynamics:AGagg} has a unique Carathéodory solution, which belongs to $S$ for all $t\geq 0$.  The solution converges to an equilibrium $\operatorname{col}\left( \bar{x},\bar{\bs\varsigma},\bar{\bs{k}},\bar{\bs z},\bar{\bs \lambda}\right) $, with  {$\psibs(\bar{x})+\bar{\bs\varsigma}=\1_N\otimes \psi (\bar{x})$}, $\bar{\bs \lambda}=\1_{N} \otimes \lambda^{*}$, and $(\bar{x},\lambda^*)$ satisfies the \gls{KKT} conditions in \eqref{eq:KKT}, hence $\bar{x}$ is the  \gls{v-GNE} of the game in \eqref{eq:game}.
	{\hfill $\square$}
\end{thm}
\begin{thm}[Convergence of Algorithm~\ref{algo:3}]\label{th:main3}
	Let Assumptions~\ref{Ass:Convexity}, \ref{Ass:StrMon}(i), \ref{Ass:Graph}, \ref{Ass:LipschitzExtPseudoagg} hold, and let $c>\underline{c}$, with $c>\underline{c}$ as in Algorithm~\ref{algo:3}.
	Then, for any initial condition in {$S={\Omega}\times \bs{E}^\perp_\nbar \times \bs{E}^\perp_{m} \times \R^{mN}_{\geq 0}$} the system in \eqref{eq:dynamics:CGagg} has a unique Carathéodory solution, which belongs to $S$ for all $t\geq 0$.  The solution converges to an equilibrium $\operatorname{col}( \bar{x},\bar{\bs\varsigma},\bar{\bs z},\bar{\bs \lambda}) $, with  {$\psibs(\bar{x})+\bar{\bs\varsigma}=\1_N\otimes\psi(\bar{x})$}, $\bar{\bs \lambda}=\1_{N} \otimes \lambda^{*}$, and  $(\bar{x},\lambda^*)$ satisfies the \gls{KKT} conditions in \eqref{eq:KKT}, hence $\bar{x}$ is the  \gls{v-GNE} of the game in \eqref{eq:game}.
	{\hfill $\square$}
\end{thm}

 \section{Multi-integrator agents}  
  \label{sec:multiintegrators}
In this section, we consider the game in \eqref{eq:game} under the following additional assumption, which is standard for \gls{NE} problems with  higher-order dynamical agents (\citealp[Asm.~1]{RomanoPavel2019}; \citealp[Def. 1]{Deng_Liang_EL_NE_2019}).
\begin{asm}\label{Ass:Unboundedfeasibleset}
	$\Omega=\R^n$.
{\hfill $\square$} \end{asm}

{Besides, we study problems where each agent is represented by a system of (mixed-order) multi-integrators:
\begin{equation}\label{eq:basemultiintegrators}
\forall \I: \qquad \left\{ x_{i,k}^{(r_{i,k})}=u_{i,k},  \quad k\in  \{1,\dots, n_i\},  \right.
\end{equation}
 where $r_{i,k}\geq 1$ and we denote by $x_{i,k}:=[x_{i}]_k$, $u_{i,k}:=[u_{i}]_k$ the  $k$-th scalar component of agent $i$ strategy and control input, respectively. }
 Our aim is to drive the agents' actions (i.e., the $x_i$ coordinates of each agent state) to a \gls{v-GNE} of the game in \eqref{eq:game}.
We emphasize that the agents are not able to directly control their strategy $x_i$ in \eqref{eq:basemultiintegrators}.
\begin{rem1}\label{rem:FL}
	{
We consider the  general form in \eqref{eq:basemultiintegrators} -- instead of  homogeneous multi-integrator systems $x_i^{(r_i)}=u_i$ as in  \cite{RomanoPavel2019} -- because
these dynamics often arise from feedback linearization of multi-input multi-output (nonlinear) systems.As an example, the feedback linearized model of a  quadrotor in \citet[Eq.~18]{LotufoColangeloNovara_TCST2020} is a combination of triple and double integrators. In general, consider any input-affine system
\begin{equation} \label{eq:nonlinear}
\forall \I: \quad \dot{z_i}=\mathfrak{f}_i(z_i)+\mathfrak{g}_i(z_i)\bar{u}_i,\quad x_i=\mathfrak{h}_i(z_i),
\end{equation}
 for smooth mappings $\mathfrak{f_i}:\R^{q_i}\rightarrow \R^{q_i}$, $\mathfrak{g_i}:\R^{q_i}\rightarrow \R^{q_i\times n_i}$, $\mathfrak{h}:\R^{q_i}\rightarrow\R^{n_i}$; the objective is to drive the controlled outputs $x_i$ to a \gls{v-GNE}. Assume that the systems in \eqref{eq:nonlinear} have, for all $z_i\in\R^{q_i}$, vector relative degree \cite[§5.1]{Isidori}  $\{r_{i,1}, \dots, r_{i,n_i}\}$, with $r_{i,1}, \dots, r_{i,n_i}\geq 1$ and $r_1+\dots +r_{n_i}=q_i$. This class includes, e.g., the Euler--Lagrangian dynamics considered in \cite{Deng_Liang_EL_NE_2019}. Then, for all $\, \I$, there is a change of coordinates $\xi_i=T_i(z_i)$ and a state feedback $\bar{u}_i=\alpha(\xi_i)+\beta(\xi_i)u_i$ such that the closed-loop system, in the new coordinates and with transformed input $u_i$, is exactly \eqref{eq:basemultiintegrators} \cite[§5.2]{Isidori}. In practice, the problem of driving the  systems in \eqref{eq:nonlinear} to a \gls{v-GNE} can be recast,  via a linearizing feedback, as that of controlling the multi-integrator agents in \eqref{eq:basemultiintegrators} to a \gls{v-GNE}.
 \hfill $\square$}
\end{rem1}
Let $\mathcal{K}_i:=\{1,\dots,n_i\}$ and 
$\mathcal{M}_i:=\{k\in \mc{K}_i \mid r_{i,k}>1\}$, for all $\I$.
 We assume that each agent is able to measure its full state.  Similarly to \cite{RomanoPavel2019}, in \eqref{eq:basemultiintegrators}, for each $\I$, we consider the controllers 
\begin{align}\label{eq:utransl}
\forall \K: \quad u_{i,k}=\tilde{u}_{i,k} - \textstyle \sum_{j=1}^{r_{i,k}-1} c_{i,k,j-1}x_{i,k}^{(j)},
\end{align}
where  $\tilde{u}_{i,k}$ is a translated input to be chosen, and $(c_{i,k,0}:=1, \dots , c_{i,k,r_{i,k}-2}, c_{i,k,r_{i,k}-1}:=1)$ are the ascending coefficients of any Hurwitz polynomial of order $(r_{i,k}-1)$, for all $i \in \mc{K}_i$.
Moreover, for all $\I$, we define the coordinate transformation
\begin{equation}\label{eq:transformation}
\col(( \col(x_{i,k},\dots,x_{i,k}^{(r_{i,k}-1)}))_{k\in\mc{K}_i})\rightarrow \col(\zeta_{i},v_{i}), 
\end{equation}
\sloppy where $v_i:=\col((v_{i,k})_{k\in\mc{M}_i})$ and $\zeta_i:=\col((\zeta_{i,k})_{k\in\mc{K}_i})$, with $v_{i,k}:=\col(x_{i,k}^{(1)}, \dots,x_{i,k}^{(r_{i,k}-1)})$, and 
\begin{equation}\label{eq:zetak}
\zeta_{i,k}:=
\left\{
\begin{aligned}
&x_{i,k}+ \textstyle \sum_{j=1}^{r_{i,k}-1} c_{i,k,j}x_{i,k}^{(j)} &&  \text{if } k\in \mc{M}_i \\
&x_{i,k} && \text{if }k\notin \mc{M}_i. 
\end{aligned}
\right.
\end{equation}
We note that, for the closed loop systems in the new coordinates, it holds, for all $i\in\I$,
%
\begin{subequations}
	\label{eq:k_multiintegrators_input}
	\begin{empheq}[left={\forall k\in \mc{M}_i:\qquad\empheqlbrace }]{align}
	\label{eq:k_multiintegrators_input:a}
	\dot{\zeta}_{i,k}&=\tilde{u}_{i,k}
	\\
	\label{eq:k_multiintegrators_input:b}
	\dot{v}_{i,k}&=E_{i,k} v_{i}+G_{i,k}\tilde u_{i,k}, \qquad 
	\end{empheq}
\end{subequations}
where $
\textstyle E_{i,k}= \begin{bmatrix}{\textstyle \0_{r_{i}-2} \ } & {I_{r_{i}-2}} \\ {1 \ } & {-c_{i,k}^\top}\end{bmatrix}\hspace{-0.2em},  \; 
\textstyle G_{i,k}= \begin{bmatrix}{\0_{r_{i}-2 }} \\ {1}\end{bmatrix}\hspace{-0.2em}, 
$
and $c_{i,k}:=\col(c_{i,k,1},\dots, c_{i,k,r_{i,k}-2})$. 

We conclude that the system in \eqref{eq:basemultiintegrators}, with the control inputs \eqref{eq:utransl}, in the new coordinates \eqref{eq:transformation}, reads as 
\begin{subequations}
	\label{eq:multiintegrators_input}
\begin{empheq}[left={\forall \I:\qquad\empheqlbrace }]{align}
\label{eq:multiintegrators_input:a}
\dot{\zeta}_{i}&=\tilde{u}_{i}
\\
\label{eq:multiintegrators_input:b}
\dot{v}_{i}&=E_{i} v_{i}+G_{i}\tilde u_{i}, \qquad 
\end{empheq}
\end{subequations}
where $\tilde{u}_i:=\col((\tilde{u}_{i,k})_{k\in\mc{K}_i})$, $E_i:=\diag((E_{i,k})_{k\in\mc{M}_i})$, $G_i: =\diag((G_{i,k})_{k\in\mc{M}_i})$, for all $\I$.

The dynamics of $\zeta_i$ in \eqref{eq:multiintegrators_input:a} are identical to the single-integrator in \eqref{eq:integrators}, with translated input $\tilde{u}_i$. 
As such, we are in a position to design $\tilde{u}_{i}$  according to Algorithm \ref{algo:2} (or \ref{algo:1}, or  \ref{algo:3} or \ref{algo:4}  for aggregative games), to drive  the variable $\zeta:=\col((\zeta_i)_{i\in\mc{I}})$ to an equilibrium $\bar{\zeta}=x^*$, where $x^*$ is the \gls{v-GNE} for the game in \eqref{eq:game}. In the following, we  show that this choice is sufficient to also control the original variables $x_i$ to the  \gls{v-GNE}.

The resulting dynamics are shown in  Algorithm~\ref{algo:mi-adaptivegain}. 
Here, $\bs{ \zeta}^i:=(\operatorname{col}(\bs{\zeta}_{j}^{i})_{j \in \mc I})$, and $\bs{\zeta}_{j}^{i}$ represents agent $i$'s estimation of the quantity
$\zeta_j$,
for $j\neq i$, while 
$\bs{\zeta}_i^i:=\zeta_i$, $\bs{\zeta}_{-i}^i:=\col(( \bs{\zeta}^i_j)_{j\in\mc{I}\backslash \{i\}})$.  Let also $\bs{\zeta}:=\col((\bs{\zeta}^i)_{\I})$. 
\begin{algorithm}[t] \caption{Multi-integrator agents (adaptive gains)}\label{algo:mi-adaptivegain}
	\vspace{0.2em}
	\textbf{Initialization}: $\forall \I$, set $\gamma_i>0$, $\bs{\zeta}^i_{-i}(0)\in \R^{n-n_i}$,  $k_i(0)\in \R$, {$z_i(0)=\0_m$}, $\lambda_{i}(0)\in \R^m$;
	\vspace{0.5em}
	\\
	\textbf{Dynamics}: $\forall \I$, for all $k\in\mc{K}_i$,
	\vspace{-0.8em}
	\vspace{0.2em}
	\begin{align*}
	\begin{aligned}
	 {x_{i,k}^{(r_{i,k})}} &= {u_{i,k}=\tilde{u}_{i,k} - \textstyle \sum_{j=1}^{r_{i,k}-1} c_{i,k,j-1}x_{i,k}^{(j)}}
	\\
	\tilde{u}_i&= \begin{multlined}[t]- \grad{x_i} J_{i}(\bs{\zeta}^{i}_{i}, \bs{\zeta}^{i}_{-i})-\diff{x_i}g_i(\bs{\zeta}_i^i)^{ \top} \lambda_{i}\\
	-
	{\textstyle \sum}_{j \in \mathcal{N}_{i}} (k_i{\rho}^i_i-k_j{\rho}^j_i)  
	\end{multlined}
	\\
	\bs{\dot \zeta}^i_{-i}&=- \textstyle\sum_{j \in \mathcal{N}_{i}} (k_i{\rho}^i_{-i}-k_j{\rho}^j_{-i})
	\\
	  \bs{\zeta}^{i}_i &=\zeta_i=\col((\zeta_{i,k})_{k\in\mc{K}_i})
	\\
	\dot k_i&=\gamma_i \|{\rho}^i \|^2 \qquad \qquad{\rho}^i=	\textstyle \sum_{j \in \mathcal{N}_{i}} (\bs{\zeta}^{j}-\bs{\zeta}^{i}) 
	\\
	{\dot z}_i&=\textstyle \sum_{j \in \mathcal{N}_{i}} \left(\lambda_{i}-\lambda_{j}\right) 
	\\
	\dot\lambda_i&=\Pi_{\R_{\geq 0}^m}(\lambda_i, g_i(\bs{\zeta}^{i}_{i})-z_i 
	-{\textstyle \sum_{j \in \mathcal{N}_{i}} } (\lambda_{i}-\lambda_{j}) )
	\\
	\end{aligned}
	\end{align*}
	\vspace{-1.1em}
\end{algorithm}
\begin{thm}[Convergence of Algorithm~\ref{algo:mi-adaptivegain}]\label{th:multiint}
	%
	%
 Let Assumptions~\ref{Ass:Convexity}, \ref{Ass:StrMon}, \ref{Ass:Graph}, \ref{Ass:Unboundedfeasibleset} hold.
 For any initial condition, the system in Algorithm~\ref{algo:mi-adaptivegain} has a unique Carathéodory solution.
The solution converges to an equilibrium $\operatorname{col}(\bar{x}, \bar{\bs{\zeta}}, \bar{\bs k}, \bar{\bs z},\bar{\bs \lambda} ) $, with  $\bar x=x^*$, $\bar{\bs{\zeta}}=\1_{N} \otimes x^{*}$, $\bar{\bs \lambda}=\1_{N} \otimes \lambda^{*}$, and   $(x^*,\lambda^*)$ satisfies the KKT conditions in \eqref{eq:KKT}, hence $x^*$ is the \gls{v-GNE} of the game in \eqref{eq:game}. 
	{\hfill $\square$}
\end{thm}

We emphasize that the proof of Theorem~\ref{th:multiint} is not based on the specific structure of Algorithm~\ref{algo:2}; in fact, the result still holds if another secondary controller with analogous convergence properties is employed to design $\tilde{u}_i$ in \eqref{eq:multiintegrators_input}.
%
%
For instance, by choosing the controller in \citet[Eq.~11]{DengNian2019}, we can address aggregative games played by multi-integrator agents over balanced digraphs. 
 \cite{RomanoPavel2019} follow a similar approach (for homogeneous multi-integrators and  \gls{NE} problems), and handle the presence of deterministic disturbances by leveraging the \gls{ISS} properties of their selected secondary controller  \cite[Eq.~47]{GadjovPavel2018}.
We have not guaranteed this robustness for our dynamics.
However, the algorithm in \cite{RomanoPavel2019} is designed for unconstrained games.
On the contrary, Algorithm~\ref{algo:mi-adaptivegain} drives the system in \eqref{eq:basemultiintegrators} to the \gls{v-GNE} of a generalized game, and ensures asymptotic satisfaction of the coupling constraints. 
We finally remark that  we  assumed
the absence of local constraints (Assumption~\ref{Ass:Unboundedfeasibleset});
nevertheless, if some are present, they can be dualized and satisfied asymptotically, as in Remark~\ref{rem:dualizedlocalconstraints}.

\section{Illustrative numerical examples}\label{sec:simulations}
\subsection{Optimal positioning in mobile sensor networks}
 \begin{figure}[t]
	\centering
	\includegraphics[width=0.92\columnwidth]{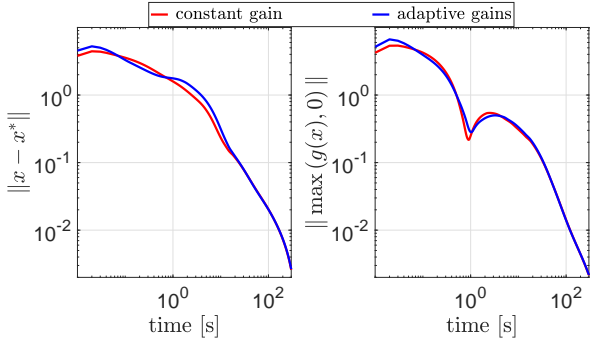}
	%
	%
	\caption{Results of Algorithms~\ref{algo:1}-\ref{algo:2} for velocity-actuated vehicles.}
	\label{fig:1}
\end{figure}
We consider a  connectivity problem formulated as a game, as in \cite{Stankovic_Johansson_Stipanovic_2012}. A group $\mc{I}=\{1,\dots,N=5\}$ of mobile sensor devices have to coordinate  their actions via wireless communication, to perform some task, e.g., exploration or surveillance. 
Mathematically, each sensor $i$ aims at autonomously  finding the position $x_i=\col(p^x_i,p^y_i)$ in a plane to optimize some private primary  objective $c_i(x_i)$, but not rolling away too much from the other devices.
%
%
This goal is represented by the following cost functions, for all $\I$:
\[ 
\textstyle J_{i}\left(x_{i}, x_{-i}\right):=c_i(x_i)+\sum_{j \in \mathcal{I}}\left\|x_{i}-x_{j}\right\|^{2}.\]
Here, $c_i(x_i):= x_{i}^{T} x_{i}+ d_{i}^\top x_i+{\sin(p^x_i)}$,
with $d_i\in\R^2$ randomly generated local parameters, for each $i\in\mc{I}$.
The useful space is restricted by the local constraints $0.1\leq p^y_i \leq 0.5$, $\forall i\in\mc{I}$. The sensors communicate over a random  undirected connected graph $\mc{G}(\mc{I},\mc{E})$. 
To preserve connectivity, the Chebyschev distance between any two neighboring agents has to be smaller than $\frac{1}{5}$, resulting in the coupling constraints  $\max\{|p^x_i-p^x_j|,|p^y_i-p^y_j|\}\leq\frac{1}{5}, \forall (i,j)\in \mc{E}$. After the deployment, all the sensors start sending the data they collect to a base station, located at  $\bar{x}=\col(0,0.3)$,  via wireless communication.
To maintain acceptable levels of transmission power consumption, the average steady state distance from the base is limited as $\textstyle \frac{1}{N} \sum_{ j \in \mathcal{N}_{i}} (x_i-\bar{x})^\top (x_i-\bar{x}) \leq \frac{1}{2}$. This setup satisfies  Assumptions~\ref{Ass:Convexity}-\ref{Ass:StrMon}. We set $c=30$ to satisfy the condition in Theorem~\ref{th:main1}; $\gamma_i=1, \forall i\in \mc{I}$; initial conditions are chosen randomly. 
We consider two different cases for the sensor physical dynamics. 

\emph{Velocity-actuated vehicles:} Each agent is a single-integrator as in \eqref{eq:integrators}. Figure~\ref{fig:1} compares the results for Algorithms~\ref{algo:1} and~\ref{algo:2} ({in logarithmic scale}) and shows convergence of both to the unique \gls{v-GNE} and asymptotic satisfaction of the coupling constraints. {In the first phase, the controllers are mostly driven by the consensual dynamics; we remark that, when the agents agree on their estimates, the two algorithms coincide}.
\begin{figure}[t]
	\centering
	\includegraphics[width=\columnwidth]{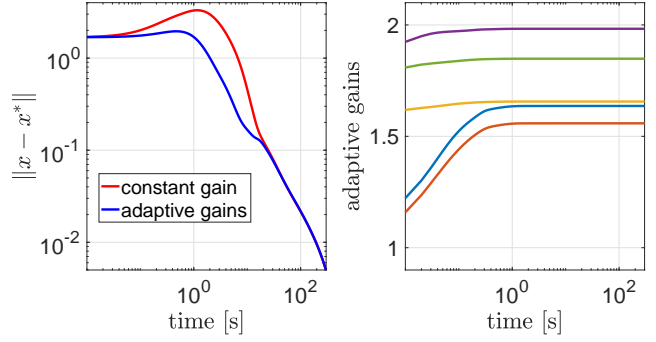}
	\caption{Results of Algorithm~\ref{algo:mi-adaptivegain} for Euler--Lagrangian vehicles linearized via feedback linearization.}
	\label{fig:2}
\end{figure}

\emph{Euler--Lagrangian vehicles}: Each  mobile sensor $\I$ is modeled as an Euler--Lagrangian systems of the form $I_i(x_i)\ddot{x}_i+C_i(x_i,\dot{x}_i)+U_i(x_i)=u_i$, where $U_i=\col(0,-1)$,
\[
\begin{aligned}
I_i(x_i) &=\left[\begin{smallmatrix} 2+0.6*\cos(p_i^y) & 0.5+0.3\cos(p_i^y)
\\  0.5+0.3\cos(p_i^y) & 0.5
\end{smallmatrix}\right], \
\\
C_i(x_i,\dot{x}_i)&=
\left[
\begin{smallmatrix}
-0.3\sin(p_{i}^y) \dot{p}_{i}^{ y} & 
-0.3\sin (p_{i }^{y})(\dot{p}_{i}^{ x}+\dot{p}_{i}^ y)
\\ 0.3\sin (p_{i }^y)  \dot{p}_{i}^x &0
\end{smallmatrix}
\right].
\end{aligned}
\]	
The systems satisfy the conditions in Remark~\ref{rem:FL} with uniform vector relative degree $\{2,2\}$. Therefore, we first apply a linearizing feedback; the problem then reduces to the control of double-integrator agents, and  we choose the transformed input (see Remark~\ref{rem:FL}) according to Algorithm~\ref{algo:mi-adaptivegain} and the analogous algorithm with constant gain (obtained by choosing $\tilde{u}_i$ in \eqref{eq:multiintegrators_input:a} according to Algorithm~\ref{algo:1}). 
The local constraints are dualized as in Remark~\ref{rem:dualizedlocalconstraints}.  
The results are illustrated in Figure~\ref{fig:2}. 
Finally, in Figure~\ref{fig:3}, we compare the trajectories of the vehicles in the velocity-actuated and Euler--Lagrangian cases.  Importantly, the local constraints are satisfied along the whole trajectory for single-integrator agents, only asymptotically for the higher-order agents.
   \begin{figure}[t]
 	\centering
 	\includegraphics[width=\columnwidth]{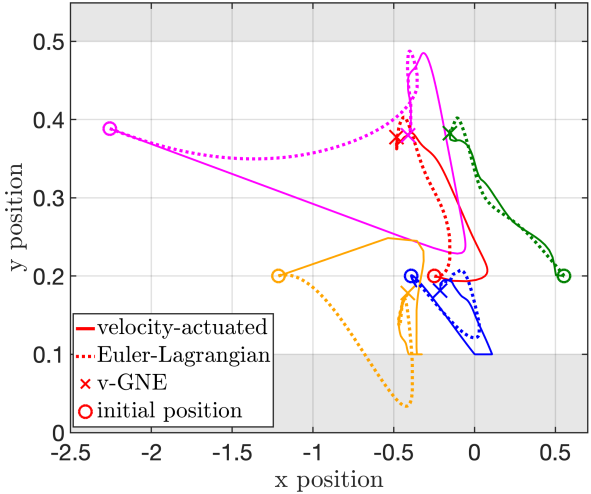}
 	\caption{Cartesian trajectories of velocity-actuated and Euler--Lagrangian vehicles, with adaptive gains.}
 	\label{fig:3}
 \end{figure}
  \subsection{Competition in power markets as aggregative game}
    \begin{figure}[t]
  	\centering
  	\includegraphics[width=\columnwidth]{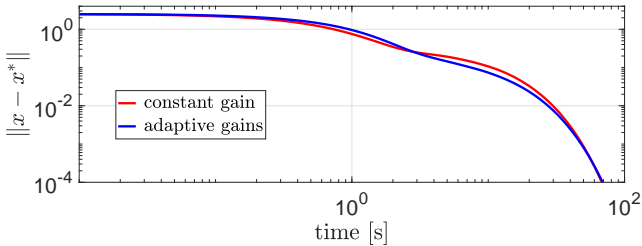}
  	\caption{{Distance from the v-GNE, for the power production in the electricity market.}}
  	\label{fig:4}
  \end{figure}  
  We consider a Cournot competition model \citep{Hoobs2007,Pavel2018}. A group  $\mc{I}=\{1,\dots,N\}$ of firms produces energy for a set of  markets $\mc{J}=\{1,\dots,m\}$, each corresponding to a different location. Each firm $\I$ controls a production plant in $n_i\leq m$ of the locations, and decides on the power outputs $x_i\in\R^{n_i}$ of its generators.
   Power is only dispatched in the location of production. Each plant has a maximal capacity, described by the local constraints $\0_{n_i} \leq x_i \leq X_i$. 
  Moreover, an \gls{ISO}  imposes an upper bound on the market share of the  producers, so that $\1_{n_i}^\top x_i\leq C_i$.
  Market clearance is guaranteed by the \gls{ISO} via external control mechanisms, but the overall power generation is bounded by markets capacities $r=\col((r_j)_{j\in\mc{J}})$. Thus, the firms share the constraints $Ax\leq r$. Here, $A=[A_1 \dots A_N]$, and $A_i\in\R^{m\times n_i}$  with
  $[A_i]_{j,k}=1$ 
  if $[x_i]_k$ is the power generated in location $j\in\mc{J}$ by agent $i$,  $[A_i]_{j,k}=0$ otherwise, for all $j\in\mc{J}$, $k=1,\dots, n_i$. In simple terms, $Ax\in \R^m$ is the vector of total power generations for each market. Each firm $\I$ aims at maximizing its profit, i.e.,  minimizing the cost 
  \[
  J_i(x_i,x_{-i})=c_i(x_i)-p(Ax)^\top A_i x_i+w(\1_{n_i}^\top  x_i),
  \]
  where $c_i(x_i)=\textstyle \sum_{k=1}^{n_i} Q_{i,k} ([x_i]_k)^2+q_{i,k}[x_i]_k$ 
  is the generation cost,  $p(Ax)^\top A_ix_i$ is the revenue, where $p:\R^m\rightarrow \R^m$ associates to each market a unitary price depending on the offer and $[p(Ax)]_j=P_j-\chi_j[Ax]_j$,
 $w(y)=w_2 y - w_1 y^2$ is a price charged by the \gls{ISO} for the use of the infrastructure.
 We set $N=20$, $m=7$ and randomly select which firms participate in each market. We choose with uniform distribution $X_i$  in $[0.3,1.3]$, $C_i$ in $[1,2]$, $r_j$ in $[1,2]$, $Q_{i,k}$ in $[8,16]$, $q_i$ in $[1,2]$, ${P}_j$ in $[10,20]$, $\chi_j$ in $[1,3]$, for all $\I$, $j\in\mc{J}$, $k=1,\dots,n_i$, $w_1$ in $[0.5,1]$, $w_2$ in $[0,0.1]$. 
 The firms cannot access the productions of all the competitors, but  can communicate with some neighbors on a connected graph. The turbine of  generator $i$ is governed by the  dynamics \citep{Deng_Liang_EL_NE_2019}
\[
\begin{aligned}
 \dot{P}_{i,k} &=\textstyle -\alpha^1_{i,k} P_{i,k}+\alpha^2_{i,k} R_{i,k} \\
\dot{R}_{i,k} &=-\alpha^3_{i,k} R_{i,k}+\alpha^4_{i,k}  u_{i,k},
\end{aligned}
\]
with $P_{i,k}=[x_i]_k$; $R_{i,k}$ and $u_{i,k}$ are the steam valve opening and  control input; the parameters  $\alpha^{\cdot}_{i,k}$'s are set as in \cite{Deng_Liang_EL_NE_2019}.
Via feedback linearization, the  problem for each generator reduces to the control of a double-integrator. 
The competition among the firms is described as an aggregative game with aggregation value  $\psi(x)=Ax$ (this is advantageous with respect to the formulation in \cite{Pavel2018}, as the firms only keep an estimate of the aggregation and firm $i$ does not need to know the quantities $A_j$, $j\neq i$). We numerically check that this setup satisfies  Assumptions \ref{Ass:StrMon}, \ref{Ass:LipschitzExtPseudoagg}. We simulate the equivalent of Algorithm~\ref{algo:mi-adaptivegain} for aggregative games, obtained by choosing $\tilde{u}_i$ in \eqref{eq:multiintegrators_input} according to Algorithms~\ref{algo:3}, \ref{algo:4}; we deal with the local constraints as in Remark~\ref{rem:dualizedlocalconstraints}. The results are shown in Figure~\ref{fig:4} and indicate fast convergence of the firms' production to the unique \gls{v-GNE}.
  
 \section{Conclusion and outlook}\label{sec:conclusion}
 
Generalized games played by nonlinear systems with maximal relative degree can be solved via continuous-time fully distributed primal-dual pseudogradient controllers, provided that the game mapping is strongly monotone and Lipschitz continuous. Convergence can be ensured even  without a-priori knowledge on the game parameters, via integral consensus. 
Seeking an equilibrium when the agents are characterized by constrained dynamics is currently an open problem. The extension of our results to the case of {direct communication}, noise and parameter uncertainties is left as future work. 
 
%

\appendix

\section{Proof of Lemma \ref{lem:equilibria-vGNE-adaptive}}\label{app:lem:equilibria-vGNE-adaptive}
 	Under  Assumption~\ref{Ass:Graph}, we have, for any $q>0$,
 	\begin{alignat}{3}
 	 \label{eq:range_Lq}
 	\operatorname{Range}\left(\L{q}\right)&= \operatorname{N u l l}\left(\1_{N}^\top \otimes I_{q}\right) &&=\bs{E}_{q}^{\perp} ,
 	\\
 	\label{eq:null_Lq}
 	\operatorname{N u l l}\left(\L{q}\right)&= \operatorname{Range}\left(\1_{N} \otimes I_{q}\right)&&=\bs{E}_{q}.
 	\end{alignat}
\indent \textit{i)} For any equilibrium $\operatorname{col}({\bar{\bs{x}}},\bs{\bar k},{\bar {\bs z}},{\bs{\bar \lambda}}) $ of \eqref{eq:dynamics:AG}, with $\zbar\in\bs{E}^\perp_m$, it holds that
\allowdisplaybreaks
	\begin{align}
\label{eq:line1}
\0  & \in \mathcal{R}^\top \! ( \bs{F} (\bar{\bs{x}} )+\mathbb{G}  (\mc{R} \bs{\bar{x}}) ^{\top} \bs{\bar \lambda}+\mathrm{N}_{ \Omega}(\mathcal{R} \bar{\bs{x}} )) \!+\!
 \bs{L}_n 
\bar K \bs{\bar{\rho}}
\\
\label{eq:line2}
\0 &=D\left(\bar {\bs{\rho}} \right)^\top \left( \Gamma \otimes I_n\right)\bar {\bs{\rho}}, \quad {\bar{\bs{\rho}}}=\bs{L}_n\bar{\bs{x}}
\\
\label{eq:line3}
\0 &=\bs{L}_{m} \bs{\bar \lambda} 
\\ 
\label{eq:line4}
\0 & \in -\bs{g}(\mc{R}\bs{\bar{x}}) + \bs{\bar{z}}  +\L{m}\bs{\bar \lambda}  +\mathrm{N}_{\R_{\geq 0}^{N m}}\left(\bs{\bar \lambda}\right),\end{align}
where $\bar K =\operatorname{diag}((\bar k_{i} I_{n})_{\I})$.
By \eqref{eq:line2} we have  $\bar{\bs{\rho}}=\0_{Nn}$, i.e., $\bar{\bs{x}}\in \bs{E}_{n}$ by \eqref{eq:null_Lq}, and by \eqref{eq:line3} and \eqref{eq:null_Lq}, we have $\bar{\bs \lambda}\in \bs{E}_{m}$. Hence,  $\bar{\bs{x}}=\1_{N} \otimes x^*$ and $\bs{\bar \lambda}=\1_{N} \otimes \lambda^*$, for some $x^*\in \R^n$, $\lambda^*\in\R^m$. 
By  left multiplying both sides of \eqref{eq:line1} by $(\1_N^\top \otimes I_n)$, by \eqref{eq:null_Lq} and since $(\1_N^\top \otimes I_n)\mc{R}^\top=I_n$, 
$\bs{F}({\1_N \otimes x^*})=F(x^*)$, $\mathcal{R} \bar{\bs{x}}=x^*$, and $\mathbb{G}  (\mc{R} \bs{\bar{x}})^\top(\1_{N} \otimes \lambda^*)=\diff{x}g(x^*)^\top \lambda^*	$,
we retrieve the first \gls{KKT} condition in \eqref{eq:KKT}.
We obtain the second condition in \eqref{eq:KKT} by   left multiplying both sides of  \eqref{eq:line4} by 
$(\1_N^\top \otimes I_m)$
and by using that 
$(\1_N^\top \otimes I_m)\bs{g}(\mc{R}\bs{\bar{x}}) =g(x^*)$, 
$(\1_N^\top \otimes I_m)\zbar=0$ and
$\textstyle {(\1_N^\top \otimes I_m)\mathrm{N}_{\R_{\geq 0}^{Nm}}(\1_N\otimes\lambda^*)=N\mathrm{N}_{\R_{\geq 0}^{ m}}(\lambda^*)=\mathrm{N}_{\R_{\geq 0}^{ m}}(\lambda^*)}$.

 \textit{ii)}  Let  $(x^*,\lambda^*)$ be any pair that satisfies the \gls{KKT} conditions in \eqref{eq:KKT}. By taking $\bar{\bs{x}}=\1_{N} \otimes x^{*}$, $\bar{\bs \lambda}=\1_{N} \otimes \lambda^{*}$ and any $\bar{\bs{k}}$, \eqref{eq:line1}-\eqref{eq:line3} are satisfied as above. It suffices to show that there exists $\bs{\bar z}\in \bs{E}_{m}^{\perp}$ such 
that \eqref{eq:line4} holds, i.e.,  that $(-\bs{g}(\mc{R}\bs{\bar{x}})+ \bar v)\in\bs{E}_{m}^\perp$, for some $\bar{v}\in\mathrm{N}_{\R_{\geq 0}^{N m}}(\bs{\bar \lambda})$. By \eqref{eq:KKT},  there exists $v^*\in \mathrm{N}_{\R^m_{\geq 0}}(\lambda^*)$ such that $g(x^*)= v^*$. 
Since  $\mathrm{N}_{\R_{\geq 0}^{N m}}(\1_{N} \otimes \lambda^{*})=\times_{i\in \mc I}\mathrm{N}_{\R^m_{\geq 0}}(\lambda^*)$,  it follows by properties of cones that $\operatorname{col}(v_1^*,\dots,v_N^*)\in \mathrm{N}_{\R_{\geq 0}^{N m}}(\bs{\bar \lambda})$ with $v_1^*=\dots=v_N^*=\frac{1}{N}v^*$. Therefore,
$(\1_{N}^\top \otimes I_{m})(-\bs{g}(\mc{R}\bs{\bar{x}})+ \bar v)=-g(x^*)+v^*=\bs 0_m$, or $(-\bs{g}(\mc{R}\bs{\bar{x}})+ \bar v)\in\bs{E}_{m}^\perp$.
{\hfill $\blacksquare$}
\section{Proof of Lemma \ref{lem:strongmon_adaptive}}\label{app:lem:strongmon_adaptive}
Let $\bs{y}=\1_n\otimes y$, for some $y \in \R^n$. We decompose $\bs{x}=\bs{x}^\perp+\bs{x}^\parallel$, where
$\bs{x}^\parallel:=P_n \x$, $\bs{x}^\perp:=P^\perp_n \x$.
Therefore, $\bs{x}^\parallel=\1_N\otimes \hat{x}$, for some $\hat x \in \R^n$. By \cite[Eq.~50]{Pavel2018}, 
\allowdisplaybreaks
\begin{align*}
\begin{multlined}
(\bs{x}-{\bs{y}})^{\top} \mathcal{R}^{\top}
(\bs{F}(\bs{x})-\bs{F}(\bs{y}) )
\geq  -\theta\|\hat{x}-y\|\|\bs{x}^{\perp}\| 
\\ 
 +\mu\|\hat{x}-y\|^{2} 
-\theta\|\bs{x}^{\perp}\|^{2} 
-\theta_{0}\|\bs{x}^{\perp}\|\|\hat{x}-y\|.
\end{multlined}
\end{align*}
For any $k^*> \underline{k}>0$, we have $K^*\succ 0$ and, by \eqref{eq:null_Lq}, 
$\operatorname{N u l l}\left( \L{n}K^*\L{n}\right)=\bs{E}_{n}$. Therefore it holds that
$
(\bs{x}  -\y)^{\top}\L{n}K^{*} \L{n} (\bs{x}-\y)\geq
 k^{*} \lambda_{2}(L)^{2}\| \bs{x} ^\perp  \|^{2}.
$
 By  $\left\|\hat{x}-y\right\|=\textstyle\frac{1}{\sqrt{N}}\left\|\bs{x}^{ \|}-\bs{y}\right\|$, we conclude that
\[
\begin{aligned}
& \quad \ \left(\bs{x}-\bs{y}\right)^\top(\mathcal{R}^\top\left(\bs{F}(\bs{x})-\bs{F}\left(\bs{y}\right)\right)+ \L{n} K^* \L{n}  \left(\bs{x}-\bs{y}\right))
\\
& \geq
\col(\|\bs{x}^\perp \|,\|\bs{x}^{ \|}-\bs{y}\| )^\top M_1 \col(\|\bs{x}^\perp \|,\|\bs{x}^{ \|}-\bs{y}\| ),
\end{aligned}
\]
with ${M}_1$ as in \eqref{eq:M2} and, for $k^*>\underline{k}$, ${M}_1\succ 0$ by Silvester's criterion. The conclusion follows since, by orthogonality, $\|\bs{x}^{ \|}-\bs{y}\|^{2}+\|\bs{x}^{\perp}\|^{2}=\|\bs{x}-\bs{y}\|^{2}$.
{\hfill $\blacksquare$}

\section{Proof of Theorem \ref{th:main2}}\label{app:th:main2}
We first rewrite the dynamics in \eqref{eq:dynamics:AG} as
\begin{align}
\label{eq:compact_adaptive}
\dot{\bs{\omega}}=\Pi_{\Xi} ({\bs{\omega}},-\mc{A}(\bs \omega)-\mc{B}(\bs \omega)),
\end{align}
where ${\bs{\omega}}:=\operatorname{col}\left( {\bs{x}},{\bs k},{\bs z},{\bs \lambda}\right)$, $\Xi:=\bs{\Omega}\times\R^{N}\times \R^{Nm}\times \R^{mN}_{\geq 0}$,
\begin{align*}
\mathcal{A}(\bs \omega):=
\left[\begin{smallmatrix}
\mathcal{R}^\top \bs F(\bs{x})+\bs{L}_n K \L{n} \x
\\ 
-D({\bs{\rho}})^{\top}\left(\Gamma \otimes I_{n}\right) {\bs{\rho}}
\\
\bs 0_{Nm}
\\
\bs{L}_{m} \bs \lambda 
\end{smallmatrix}\right]\hspace{-0.3em}, \
\mc{B}(\bs \omega):=
\left[
\begin{smallmatrix}
 {{\mathcal{R}^\top \mathbb{G} (\mc{R}\x)^\top} \bs{\lambda}}
\\ 
\0_N
\\
-\bs{L_m} \bs{\lambda}
\\ 
{-\bs{g}(\mc{R}\x)+\bs{z}}
\end{smallmatrix} \right]\hspace{-0.3em}. 
\end{align*}
By  Assumption~\ref{Ass:Convexity} and Lemma~\ref{lem:LipschitzExtPseudo}, $\mc{A}$ and $\B$ are locally Lipschitz; therefore, for any initial condition in $\Xi$, the system \eqref{eq:compact_adaptive} has a unique \emph{local} Carathéodory  solution, contained in $\Xi$ \citep{Cherukuri2016}. Moreover, we note that the set $S=\{\bs{w}\in \Xi \mid \bs{z}\in\bs{E}^\perp_m \}$ is invariant for the system \eqref{eq:compact_adaptive}, since for all $\bs{\omega}\in\Xi$,
$\textstyle \frac{\partial }{\partial \w}((\1_N\otimes I_m)^\top \bs{z}) \dot{\bs{\omega}}= (\1_N^\top \otimes I_m) \bs{L}_m \bs{\lambda}=\0_m$.

Let  $\Phi:=P_m+\bs{L}^+_m$, where  $\bs{L}^+_m$ is the Moore-Penrose pseudo-inverse of $\bs{L}_m$, and we recall that $P_m=\textstyle \frac{1}{N}\1_N\1_N^\top \otimes I_m$ is the projection matrix on $\bs{E}_m$. By properties of the pseudo-inverse and \eqref{eq:null_Lq}, 
$\bs{L}^+_m={\bs{L}^+_m}^\top$, $\bs{L}^+_m\succeq 0$ and $\operatorname{N u l l}(\bs{L}^+_m)=\bs{E}_m$. Since $\operatorname{N u l l}(P{_m})=\bs{E}_m^\perp$ and $P_m\succeq 0$, we have  $\Phi\succ 0$. Also, $\bs{L}^+_m \bs{L}_m=I_{Nm}-P_m=P^\perp_m$ is the projector matrix on $\operatorname{Range} (\bs{L}_m)=\bs{E}^\perp_m$.  We define the quadratic Lyapunov function 
\begin{align*}
V&=\textstyle \frac{1}{2}\|\bs{\omega}-\bar{\bs\omega} \|^2_Q:=(\w-\wbar)^\top Q (\w-\wbar)
\\
 &=\textstyle \frac{1}{2}(\| \bs{x}-\bar{\bs{x}} \|^2 +\| \bs{k}-\bs{\bar k} \|^2_{\Gamma^{-1}} +\| \bs{z}-\bs{\bar z} \|_{{\Phi}}^2 +\| \bs{\lambda}-\bs{\bar \lambda} \|^2 ),
\end{align*}
where $Q=:\operatorname{diag}(I_{Nn},\Gamma^{-1},{\Phi},I_{Nm})$, and $\bar{\bs{x}}=\1_{N} \otimes x^{*}$, $\bar{\bs \lambda}=\1_{N} \otimes \lambda^{*}$, where the pair $(x^*,\lambda^*)$ satisfies the \gls{KKT} conditions in \eqref{eq:KKT}, $\bar{\bs k}$ such that $k^*:=\min(\bar{\bs{k}})\geq \underline{k}$, with $\underline{k}$ as in \eqref{eq:M2}, $\bar{\bs z}\in\bs{E}^\perp_m$ chosen such that $\bar{\bs{\omega}}:=\operatorname{col}\left( \bar{\bs{x}},\bar{\bs k},\bar{\bs z},\bar{\bs \lambda}\right)$ is an equilibrium for \eqref{eq:dynamics:AG}, and such a $\bar{\bs z}$ exists by the proof of Lemma \ref{lem:equilibria-vGNE-adaptive}. 
%
Therefore, for any $\bs{\omega}\in S$, we have
\begin{align}
\nonumber
\dot V(\bs \omega ):&=\nabla V(\bs{\omega}) \dot{\bs{\omega}}=(\bs{\omega}-\bar{\bs{\omega}})^{\top} Q\dot{\bs{\omega}}= 
\\
\nonumber
&= (\bs{\omega}-\bar{\bs{\omega}})^{\top} Q \Pi_{\Xi} ({\bs{\omega}},-\mathcal{A}(\bs \omega)-\mathcal{B}(\bs \omega))
\\
\label{eq:usefulstep}
&\leq   (\bs{\omega}-\bar{\bs{\omega}})^{\top} Q (-\mathcal{A}(\bs \omega)-\mc{B}( \bs{\omega})),
\end{align}
where the last inequality follows from Lemma~\ref{lem:minoratingproj} and by exploiting the structure of $Q$ and $\Xi$.
By Lemma~\ref{lem:minoratingproj}, it also holds that
$
(\bs{\omega}-\bar{\bs{\omega}})^{\top} Q (-\mathcal{A}(\bar{\bs \omega})-\mc{B} (\bar{\bs{\omega}}))\leq 0.
$
By subtracting this term from \eqref{eq:usefulstep}, we obtain 
$$ \dot V(\bs \omega ) \leq  -(\bs \omega - \bs{\bar \omega })^\top Q ( \mathcal{A}(\bs \omega ) -\mathcal{A}(\bs{ \bar  \omega})+\mathcal{B}(\bs \omega ) -\mathcal{B}(\bs{ \bar  \omega})  ). $$
 Besides, for any  $\bs{z}\in\bs{E}^\perp_m$,  by $\bs{L}^+_m \bs{L}_m=P^\perp_m$ and \eqref{eq:null_Lq}, we have  $(\bs{z}-\bs{\bar{z}})^\top\Phi \bs{L}_m (\bs{\lambda}-\bar{\bs{\lambda}})=(\bs{z}-\bs{\bar{z}})^\top (\bs{\lambda}-\bar{\bs{\lambda}})$, and hence \allowdisplaybreaks
\[
\begin{aligned}
& \quad \ (\bs \omega - \bs{\bar \omega })^\top Q (\mathcal{B}(\bs \omega ) -\mathcal{B}(\bs{ \bar  \omega}) ) \\
& = 
 (\x-\bar{\x})^\top  \mathcal{R}^\top (\mathbb{G} (\mc{R}\x)^\top\bs{\lambda}-
{\mathbb{G} (\mc{R}\bar{\x})^\top}\bar{\bs{\lambda}})
\\
& \quad \ +(\bs{\lambda}-\bar{\bs{\lambda}})^\top(-\bs{g}(\mc{R}\x)+\bs{g}(\mc{R}\bar{\x}))
\\&=(x-x^*)^\top(\nabla_{\!\! y}(\bs{g} (y)^\top\bs{\lambda}) |_{y=x}-\nabla_{\!\! y} (\bs{g} (y)^\top\bs{\bar{\lambda}})|_{y=x^*})
\\ & \quad \ -(\bs{\lambda}-\bs{\bar{\lambda}})^\top(\nabla_{\!\! y}(\bs{g} (x)^\top y) |_{y=\bs{\lambda}}-\nabla_{\!\! y} (\bs{g} (x^*)^\top y)|_{y=\bar{\bs{\lambda}}}) \geq 0\\
\end{aligned}
\]
and the last inequality holds, for any $\bs{\omega}\in S$, by applying \cite[Th.~1]{Rockafellar_saddle} (since $\l,\lbar \in\R_{\geq 0}^{Nm}$ and by  Assumption~\ref{Ass:Convexity}). Therefore, for any $\bs{\omega}\in S$, it holds that:
\begin{align}
\nonumber
\dot V(\bs \omega )&\leq  -(\bs \omega - \bs{\bar \omega })^\top Q\left( \mathcal{A}(\bs \omega ) -\mathcal{A}(\bs{ \bar  \omega})\right)
\\
\label{eq:usefulstep2}
&  \begin{aligned}
&=-(\bs{x}-\bar{\bs{x}})^{\top} \mathcal{R}^{\top}
\left( \bs{F}(\bs{x})-\bs{F}(\bs{\bar x}) \right) 
\\
&\quad -(\bs{x}-\bar{\bs{x}})^{\top} (\bs{L}_n K \bs{L}_n(\bs{x}-\bar{\bs{x}}))
\\
& \quad  +(\bs{k}-\bar{\bs{k}})^{\top} \Gamma^{-1} D({\bs{\rho}})^{\top}\left(\Gamma \otimes I_{n}\right) {\bs{\rho}}\\
& \quad -(\bs \lambda -\bar{\bs \lambda})^\top \bs{L}_{m} (\bs \lambda -\bar{\bs \lambda}),
\end{aligned}
\end{align}
where we used that 
 $\bar{\bs\rho}:=\bs{L}_{n} \bar{\bs{x}}=0$. 
 For the last addend in \eqref{eq:usefulstep2}, we can write
 $(\bs \lambda -\bar{\bs \lambda})^\top \bs{L}_{m} (\bs \lambda -\bar{\bs \lambda})=\bs \lambda^\top \bs{L}_{m} \bs \lambda $ by \eqref{eq:null_Lq} and,  by \cite[Th.~$18.15$]{Bauschke2017}, $\bs \lambda^\top \bs{L}_{m} \bs \lambda \geq \frac{1}{2\uplambda_{\textnormal{max}}(L)} \| \bs{L}_{m} \bs{\lambda} \|^2 $.
The third addend in \eqref{eq:usefulstep2} can be rewritten as
$
(\bs{k}-\bar{\bs{k}})^{\top} \Gamma^{-1} D({\bs{\rho}})^{\top}\left(\Gamma \otimes I_{n}\right) {\bs{\rho}}
=\textstyle  \sum_{i=1}^{N}\left(k_{i}-\bar{k}_i\right) {\bs{\rho}}^{i \top} {\bs{\rho}}^{i} 
= {\bs{\rho}}^{\top}(K-\bar{K}) {\bs{\rho}} 
= \bs{x}^{\top}\bs{L}_n(K-\bar{K}) \bs{L}_n  \bs{x}
= (\bs{x}-\bar{\bs{x}})^\top\bs{L}_n(K-\bar{K}) \bs{L}_n (\bs{x}-\bar{\bs{x}}),
$
where $\bar{K}:=\operatorname{diag}((\bar{k}_i I_n)_{\I})$. Therefore, the sum of the second and third term in \eqref{eq:usefulstep2} is 
$-(\bs{x}-\bar{\bs{x}})^{\top} \bs{L}_n \bar{K} \bs{L}_n (\bs{x}-\bar{\bs{x}})\leq-(\bs{x}-\bar{\bs{x}})^{\top}\bs{L}_n K^{*} \bs{L}_n (\bs{x}-\bar{\bs{x}})$, where $K^*:=k^*I_{Nn}$. By Lemma \ref{lem:strongmon_adaptive}, we finally get
\begin{align}
\label{eq:upperbound_adaptive}
\dot V
\leq
-  \uplambda_{\textnormal{min}}(M_1) \|  \bs{x}-\bar{\bs{x}}\|^2 
-\textstyle \frac{1}{2\uplambda_{\textnormal{max}}(L)} \| \bs{L}_{m} \bs{\lambda} \|^2,
\end{align}
with $M_{1}\succ 0$ as in Lemma \ref{lem:strongmon_adaptive}.

Let $\mathcal{\bar{P}}$ be any compact sublevel set of $V$ (notice that $V$ is radially unbounded) containing the initial condition $\bs{\omega}(0)\in S$. $\mathcal{\bar{P}}$ is invariant for the dynamics, since $\dot{V}(\bs{\omega}) \leq 0$ by \eqref{eq:upperbound_adaptive}. The set $\mc{P}:=\mathcal{\bar{P}}\cap S$ is compact, convex and invariant, therefore, by exploiting Lemma~\ref{lem:LipschitzExtPseudo} and the continuous differentiability in  Assumption~\ref{Ass:Convexity}, we conclude that $\mathcal{A}+\mc{B}$ is Lipschitz continuous on $\mc{P}$. Therefore, for any initial condition, there exists a unique \emph{global} Carathéodory solution to \eqref{eq:dynamics:AG}, that belongs to $\mathcal{P}$ (and therefore is bounded) for every $t$ \cite[Prop.~$2.2$]{Cherukuri2016}. Moreover,  by  \cite[Th.~$2$]{DePersisGrammatico2018}, the solution converges to the largest invariant set $\mc{O} \subseteq \{\bs\omega \in \mc{P}  \text{s.t.}\dot V(\bs \omega)=0\}$. 

We can already conclude  that $\bs{x}$ converges to the point $\1_N\otimes x^*$, with $x^*$  the unique \gls{v-GNE} of the game in \eqref{eq:game}. We next show convergence  of the other variables.  
Take any point $\underline{\w}:=\operatorname{col}( \underline{\bs{x}},\underline{\bs k},\underline{\bs z},\underline{\bs \lambda}) \in \mc{O}$. Since $\dot{V}(\underline{\w})=0$, by \eqref{eq:upperbound_adaptive} we have  $\underline{\bs{x}}=\xbar=\1_N\otimes x^*,$ and $\underline{\bs \lambda}\in \bs{E}_{m}$, i.e.  $\underline{\bs \lambda}=\1_N \otimes \underline\lambda$, for some  $\underline\lambda \in \R^m_{\geq 0}$. Therefore, by expanding \eqref{eq:usefulstep}, by  $\underline{\bs{x}}=\bar{\bs{x}}$,  $\underline{\bs{\rho}}:=\bs{L}_{n} \underline{\bs{x}}=\0_{Nn}$ and  \eqref{eq:null_Lq}, we have
\begin{align}
\nonumber 
0 &=(\bs{\underline\lambda}-\bs{\bar \lambda})^\top (\bs{g}(R\bar{\bs{x}})- {\bs{\underline z} } )
=
(\underline\lambda- \lambda^{* })^\top g(x^*)
\\
\label{eq:usefulstep3}
&=\underline\lambda^\top g(x^*)
=\bs{\underline\lambda}^\top (\bs{g}(\mathcal R\bar{\bs{x}}) {- \bs{\underline z}) } ),
\end{align}
where in the second equality we have used that $\underline{\bs{z}}\in \bs{E}^\perp_m$ \ and the third equality follows from the  KKT conditions in \eqref{eq:KKTb}. 
Then, let $\underline{\bs{\omega}}(t)=\operatorname{col}( \underline{\bs{x}}(t),\underline{\bs{k}}(t), \underline{\bs z}(t),\underline{\bs \lambda}(t))$ be the trajectory of \eqref{eq:compact_adaptive} starting at  $\underline{\w}$. By invariance of $\mc{O}$,  $\underline{\bs{x}}(t)=\bar{\bs{x}}$ and $\underline{\lambda}(t)\in\bs{E}_m$, for all $t$. 
 Therefore, by \eqref{eq:dynamics:AG:2}-\eqref{eq:dynamics:AG:4}, it holds that ${\underline{\bs k}}(t)\equiv{\underline{\bs k}}$,  $\underline{\z}(t)\equiv\underline{\z}$, for all $t$.  
 Hence, the quantity $v:= (\bs{g}(\mathcal R\underline{\bs{x}}(t)) {-\bs{L}_{m} \underline{\bs{\lambda}}(t)- \underline{\bs{z}}(t) } )$ is a constant along the trajectory $\underline{\bs{\omega}}(t)$.
Suppose by contradiction that $[v]_k>0$, for some $k=\{1,\dots,Nm\}$. Then,  $[\dot{\underline{\bs \lambda}}(t)]_k=[v]_k$ for almost all  $t$, by \eqref{eq:dynamics:AG:5}, and ${\underline{\bs \lambda}}(t)$ grows indefinitely. Since all the solutions of \eqref{eq:dynamics:AG} are bounded, this is a contradiction. Therefore, $v\leq \0_m$, and $\underline{\bs \lambda}^\top v=0$ by \eqref{eq:usefulstep3}. Equivalently, $\textstyle v\in {N}_{\R_{\geq 0}^{N m}}(\underline{\bs \lambda})$, hence ${\underline{\bs \lambda}}(t)\equiv \underline{\l}$, for all $t$. 
We conclude that the points in $\mc{O}$ are equilibria.
%
%
%

Moreover, the set $\Lambda (\bs{\omega}_0)$ of $\omega$-limit points\footnote{$z :[0, \infty) \rightarrow \mathbb{R}^{n}$ has an $\omega$-limit point at $\bar{z}$ if there exists a nonnegative diverging sequence $\{t_k\}_{k\in\N}$  such that $ z\left(t_{k}\right) \rightarrow \bar{z}$}of the solution to \eqref{eq:dynamics:AG} starting at any $\bs{\omega}_0\in S$ is nonempty (by Bolzano-Weierstrass theorem, since all the trajectories of \eqref{eq:dynamics:AG} are bounded), and $\Lambda(\bs{\omega}_0)\subseteq \mc{O}$ (see the proof of \cite[Th.2]{DePersisGrammatico2018}). Hence, all the $\omega$-limit points are equilibria. We  next show that, for any for any $\bs{\omega}_0\in S$, $\Lambda(\bs{\omega}_0)$ is a singleton; as a consequence, the solution converges to that point \cite[Lemma $1.14$]{Bauschke2017}.
 For the sake of contradiction, assume  that there are two $\omega$-limit points $\bs{\omega}_1= \operatorname{col}( {\bar{\bs{x}}},\hat{{\bs k}},{{\bs{z}}}_1,{{\bs \lambda}}_1) $, $\bs{\omega}_2= \operatorname{col}( {\bar{\bs{x}}},\hat{{\bs k}},{{\bs{z}}}_2,{{\bs \lambda}}_2) \in\Lambda(\bs{\omega}_0)$, with $\bs{\omega}_1\neq \bs{\omega}_2$. We note that $\bs{\omega}_1$ and $\bs{\omega}_2$ must have the same vector of adaptive gains $\hat{{\bs k}}$  by definition of $\omega$-limit point, since the $k_i$'s in Algorithm~\ref{algo:2} are nonincreasing. Let  $\bs{\omega}_3= \operatorname{col}( {\bar{\bs{x}}},\hat{{\bs k}}+\1\alpha,{{\bs{z}}}_1,{{\bs \lambda}}_1)$, $\alpha\in \R$ chosen such that $\min(\hat{{\bs k}}+\1\alpha)>\underline{k}$, $\underline{k}$ as in \eqref{eq:M2}. By \eqref{eq:upperbound_adaptive}, $\| \bs{\omega}(t)-\bs{\omega}_3\|_Q$ is nonincreasing along the trajectory $\w(t)$ of \eqref{eq:compact_adaptive} starting at $\w_0$. Thus, by definition of $\omega$-limit point, it holds that $\|\bs{\omega}_1-\bs{\omega}_3\|_Q=\|\bs{\omega}_2-\bs{\omega}_3\|_Q$, or $\|\col (\0_{Nn},\alpha\1_N,\0_{Nm},\0_{Nm})\|_Q=\|\col (\0_{Nn},\alpha\1_N,{\bs \lambda}_1-{\bs \lambda}_2,{{\bs{z}}}_2-{{\bs{z}}}_1)\|_Q$. Equivalently, $\bs{\omega}_1=\bs{\omega}_2$, that is a contradiction.
%
{\hfill $\blacksquare$}

\section{Proof of Theorem~\ref{th:main1}}\label{app:th:main:1}
The proof follows as for Theorem~\ref{th:main2}, by defining ${\bs{\omega}}:=\operatorname{col}\left( {\bs{x}},{\bs z},{\bs \lambda}\right)$, $\Xi:=\bs{\Omega}\times \R^{Nm}\times \R^{mN}_{\geq 0}$,
\begin{align*}
\mathcal{A}(\bs \omega):=
\left[\begin{smallmatrix}
\mathcal{R}^\top \bs F(\bs{x})+c\bs{L}_n \x
\\
\bs 0_{Nm}
\\
\bs{L}_{m} \bs \lambda 
\end{smallmatrix}\right]\hspace{-0.3em}, \
\mc{B}(\bs \omega):=
\left[
\begin{smallmatrix}
{\mathcal{R}^\top \mathbb{G} (\mc{R}\x)^\top} \bs{\lambda}
\\ 
-\bs{L_m} \bs{\lambda}
\\ 
-\bs{g}(\mc{R}\x)+\bs{z}
\end{smallmatrix} \right]\hspace{-0.3em},
\end{align*}
with Lyapunov function $V(\w)=\textstyle \frac{1}{2}(\| \bs{x}-\bar{\bs{x}} \|^2 +\| \bs{z}-\bs{\bar z} \|_{{\Phi}}^2 +\| \bs{\lambda}-\bs{\bar \lambda} \|^2 )$,
and by exploiting, in place of Lemma~\ref{lem:strongmon_adaptive}, Lemma~3 in \cite{Pavel2018}.
{\hfill $\blacksquare$}

\section{Proof of Lemma~\ref{lem:strongmon_adaptive_agg}}\label{app:lem:strongmon_adaptive_agg}
By   Assumptions~\ref{Ass:StrMon} and ~\ref{Ass:LipschitzExtPseudoagg}, we have
\[
\begin{aligned}
& \quad \ (x-x')^\top (\bs{\tilde{F}}(x,\bs{\sigma})-\bs{\tilde{F}}(x^\prime,\sigmabs^\prime)
\\
&=
\begin{multlined}[t]
 (x-x')^\top (\bs{\tilde{F}}(x,\bs{\sigma})-\bs{\tilde{F}}(x,P_\nbar \psibs(x)) 
\\
 ~~~~~~~~~~~~+ \bs{\tilde{F}}(x,P_\nbar \psibs(x))-\bs{\tilde{F}}(x^\prime, P_\nbar \psibs(x^\prime )) )
\end{multlined}
\\
 &\geq  {\mu} \|x-x^\prime\|^2-\tilde{\theta}_\sigma\|x-x^\prime\| \|{\bs\sigma}-\1_N\otimes \psi(x)\|.
\end{aligned}
\]
 Moreover, we note that  $(\bs{\sigma}-\1_N\otimes\psi (x))\in \bs{E}^{\perp}_{\bar{n}}$, since  $P_\nbar\bs{\sigma}=P_\nbar \psibs (x)$, and  ${\bs{\sigma}}^\prime\in\bs{E}_{\bar{n}}$. Hence, by \eqref{eq:null_Lq}, we have  
$
(\bs {\sigma}-{\bs{\sigma}}^\prime)^\top\L{\nbar}K^*\L{\nbar}(\bs {\sigma}-{\bs{\sigma}}^\prime) \geq 
k^*\uplambda_2(L)^2\| \bs{\sigma}-P_\nbar \psibs(x)\|^2,
$
and the conclusion follows readily.
{\hfill $\blacksquare$}

\section{Proof of Theorem \ref{th:main4} }\label{app:th:main4}
The dynamics  \eqref{eq:dynamics:AGagg} can be recast in the form \eqref{eq:compact_adaptive}, with  ${\bs{\omega}}=\operatorname{col}\left( x,\bs\varsigma,{\bs k},{\bs z},{\bs \lambda}\right)$, $\Xi=\Omega\times \R^{N\nbar}\times \R^N\times \R^{Nm}\times\R_{\geq 0}^{Nm}$, 
\begin{align*}
\mathcal{A}(\bs \omega)=\hspace{-0.2em}
\left[
\begin{smallmatrix}
\bs {\tilde{F}}(x,\bs{\sigma})+{B^\top} \L{\nbar} K \L{\nbar} \sigmabs
\\
\L{\nbar} K \L{\nbar} \sigmabs 
\\ 
-D({\bs{\rho}})^{\top}\left(\Gamma \otimes I_{\bar{n}}\right) {\bs{\rho}}
\\
\bs 0_{Nm}
\\
\bs{L}_{m} \bs \lambda
\end{smallmatrix}
\right]\hspace{-0.3em}, \;
\B(\w)
=\hspace{-0.2em}
\left[
\begin{smallmatrix}
\mathbb{G}(x)^\top \l
\\ 
\0_{N\nbar}
\\ 
\0_N
\\
-\L{m}\l
\\ 
-\bs{g}(x)+\z
\end{smallmatrix}
\right]\hspace{-0.3em}.
\end{align*}
By proceeding as in the proof of Theorem~\ref{th:main2},
we note that the set $S$ is invariant for the dynamics, since, for all $\w\in S$,
 $\textstyle \frac{\partial}{\partial \w}(P_{\nbar}\varsigmabs) \dot{ {\w}}=\0_{N\nbar}$,
 $\textstyle \frac{\partial}{\partial \w}(P_{m}\z) \dot{ {\w}}=\0_{Nm}$.

 Analogously to the proof of Lemma~\ref{lem:equilibria-vGNE-adaptive}, it can be shown that any equilibrium point ${\bar{\bs{\omega}}}:=\operatorname{col}\left( \bar{x},\bar{\bs\varsigma},{\bar{\bs k}},{\bar{\bs z}},{\bar{\bs \lambda}}\right)\in S$ 
of \eqref{eq:dynamics:AGagg}  is such that $\bar{\bs{\lambda}}=\1_N\otimes \lambda^*$, the pair $(\bar{x},\lambda^*)$ satisfies the KKT conditions in \eqref{eq:KKT}, and
 {$\bar{{\bs\sigma}}:=\psibs(\bar{x})+\bar{\bs\varsigma}=\1_N\otimes \psi(\bar{x})$}. Moreover, for any pair $(x^*,\lambda^*)$ satisfying the KKT conditions in \eqref{eq:KKT}, there exists $\bar{\bs{z}}\in \R^{mN}$ such that $\operatorname{col}(x^*,{\1_N\otimes \psi({x}^*)-\psibs(x^*)},\bs{k},\bar{\bs{z}},\1_N\otimes \lambda^*)\in S$ is an equilibrium for  \eqref{eq:dynamics:AGagg}, for any $\bs{k}\in\R^N$. 
 The proof is omitted because of space limitations. 

 Let ${\bar{\bs{\omega}}}=\operatorname{col}\left( \bar{x},\bar{\bs\varsigma},{\bar{\bs k}},{\bar{\bs z}},{\bar{\bs \lambda}}\right)\in S$ be an equilibrium of \eqref{eq:dynamics:AGagg} such that $k^*=\min(\bar{\bs{k}})>\underline{k}$, $\underline{k}$ as in \eqref{eq:M4}, and consider the quadratic Lyapunov function $V=\textstyle \frac{1}{2}\|\bs{\omega}-\bar{\bs\omega} \|^2_Q$,
 where $Q=\operatorname{diag}(I_{n},I_{N\nbar},\Gamma^{-1},P_m+\L{m}^+,I_{Nm})$.
Analogously to the proof of Theorem~\ref{th:main1}, it holds that $(\bs \omega - \bs{\bar \omega })^\top Q (\mathcal{B}(\bs \omega ) -\mathcal{B}(\bs{ \bar  \omega}))\geq 0$, and that $\dot V(\bs \omega )\leq  -(\bs \omega - \bs{\bar \omega })^\top Q\left( \mathcal{A}(\bs \omega ) -\mathcal{A}(\bs{ \bar  \omega})\right)$, for all $\w\in S$.
Also we note that
\[
\begin{aligned}
& (x-\bar{x})^\top B^\top \L{\nbar}K\L{\nbar}(\sigmabs -\bar{\sigmabs})+(\bs\varsigma-\bar{\bs\varsigma})^\top \L{\nbar}K\L{\nbar}(\sigmabs -\bar{\sigmabs})\\
 &=(\bs\varsigma+Bx+d-(\bar{\bs\varsigma}+B\bar{x}+d)) ^\top\L{\nbar}K\L{\nbar}(\sigmabs -\bar{\sigmabs})\\
 &=(\sigmabs -\bar{\sigmabs})^\top\L{\nbar}K\L{\nbar}(\sigmabs -\bar{\sigmabs}),
\end{aligned}
\]
where $d:=\col((d_i)_{\I}))$, and that $(\bs{k}-\bar{\bs{k}})^{\top} \Gamma^{-1} D({\bs{\rho}})^{\top}(\Gamma \otimes I_{n}) {\bs{\rho}}
= (\bs{\sigma}-\bar{\bs{\sigma}})^{\top}\bs{L}_n(K-\bar{K}) \bs{L}_n (\bs{\sigma}-\bar{\bs{\sigma}})$ as in the proof of Theorem~\ref{th:main1}. Hence, by Lemma~\ref{lem:strongmon_adaptive_agg}, we obtain, for all $\w\in S$ 
\begin{align*}
\begin{multlined}
\dot V(\bs \omega) 
\leq- \uplambda_{\textnormal{min}}(M_2) (\|  x-{\bar x}\|^2+  \|  \bs\sigma-\1_N \otimes \psi(x) \|^2)
\\-
\textstyle\frac{1}{2\uplambda_{\textnormal{max}}(L)}\| \bs{L}_{m} \bs \lambda \|^2 ,
\end{multlined}
\end{align*} 
 with $M_2\succ 0$ as in \eqref{eq:M4}. 
Then,  existence of a unique global  solution for the  system in \eqref{eq:dynamics:AGagg} and  convergence to an equilibrium point
 follows as for Theorem \ref{th:main2}.
{\hfill $\blacksquare$}

\section{Proof of Theorem \ref{th:main3} }\label{app:th:main3}
Analogously to Lemma~\ref{lem:strongmon_adaptive_agg}, it can be shown that 
for any $c>\underline{c}$, for any $(x,\bs{\sigma})$ such that $P_{\nbar}\bs{\sigma}=P_{\nbar} \psibs (x)$
and any $(x^\prime,\bs{{\sigma}}^\prime)$ such that ${\bs{\sigma}}^\prime=P_{\nbar} \psibs(x^\prime)$, it 
holds  that ${(x-x^\prime)^\top(\bs{\tilde{F}}(x,{\bs\sigma})-	\bs{\tilde{F}}(x^\prime,{{\bs\sigma}^\prime})) }+{ c({\bs\sigma}-{{\bs\sigma}^\prime})^\top\L{\nbar}({\bs\sigma}-{{\bs\sigma}^\prime})} \allowbreak \geq \delta
	\|
	\col (
	x-x^\prime,
	{\bs\sigma}
	-{\1_N\otimes \psi(x) })
	\|^2$, for some $\delta>0$.
Then, the proof follows analogously to Theorem~\ref{th:main4}.
 \hfill $\blacksquare$

\section{Proof of Theorem~\ref{th:multiint}}
 	 Under the coordinate transformations in \eqref{eq:transformation}, the dynamics in Algorithm~\ref{algo:mi-adaptivegain} read as  \eqref{eq:multiintegrators_input}, where the input $\tilde{u}_i$ in Algorithm~\ref{algo:mi-adaptivegain} has been chosen by design according to Algorithm~\ref{algo:2}, under Assumption~\ref{Ass:Unboundedfeasibleset}. 
 Therefore, existence of a unique bounded solution and convergence of ${\zeta}_i$ to $x_i^*$ (and of the variables $\bs{\zeta}^i, k_i, z_i,\lambda_i$), for all $\I$, follows from Theorem~\ref{th:main2}.
 {On the other hand, we note that, for all $\I$ and all $k\in\mc{M}_i$, $E_{i,k}$  is Hurwitz, because it is in canonical controllable form and the coefficients of the last row are by design the coefficients of an Hurwitz polynomial. 
 	Therefore, $E_i$ is also Hurwitz,} and  hence the  dynamics in \eqref{eq:multiintegrators_input:b} are \gls{ISS} with respect to the input $\tilde{u}_i$ \cite[Lemma 4.6]{Khalil}. In turn, the input $\tilde u_i$ is bounded, by boundedness of trajectories in Theorem~\ref{th:main2},  Assumption~\ref{Ass:Convexity} and Lemma~\ref{lem:LipschitzExtPseudo}; moreover, by the convergence in Theorem~\ref{th:main2}, the \gls{KKT} conditions in \eqref{eq:KKT} and by 
 continuity, we have that  $\tilde u_i\rightarrow \0_{n_i}$ for $t\rightarrow \infty$.
 Hence, for all $\I$, $v_i\rightarrow \0$ asymptotically (this follows by definition of \gls{ISS}, see \cite[Ex.~$4.58$]{Khalil}). By the definition of ${\zeta}_i$, we  also have $x_i\rightarrow x_i^*$, for all $\I$. 
 \hfill$\blacksquare$


\bibliographystyle{model5-names.bst}
{\scriptsize\bibliography{library}}

\end{document}